\documentclass[11pt,reqno]{amsart}

\newcommand{\nc}{\newcommand}
\nc{\R}{\mathbb{R}}
\nc{\C}{\mathbb{C}}
\nc{\Z}{\mathbb{Z}}
\nc{\Pp}{\mathbb{P}}
\nc{\EE}{\mathbb{E}}
\nc{\cA}{\mathcal{A}}
\nc{\cC}{\mathcal{C}}

\nc{\by}{\bm{y}}
\nc{\bx}{\bm{x}}
\nc{\bn}{\bm{n}}
\nc{\br}{\bm{r}}
\nc{\bth}{{\bm{\theta}}}
\nc{\la}{\lambda}
\nc{\om}{\omega}
\nc{\obx}{\overline{\bx}}
\nc{\tbx}{\widetilde{\bx}}
\nc{\oth}{\overline{\bm{\theta}}}
\nc{\tth}{\widetilde{\bm{\theta}}}
\nc{\bvth}{\bm{\vartheta}}
\nc{\bet}{{\bm{\xi}}}
\nc{\bzet}{{\bm{\zeta}}}
\nc{\bxi}{{\bm{\xi}}}
\nc{\bys}{\by^{R}}
\nc{\ys}{{\it y}^R}
\nc{\IMj}{\mathcal{I}^{{M}}_{_j}\hspace{-0.02in}}
\nc{\ICINTj}{\mathcal{I}^{{CINT}}_{_j}\hspace{-0.02in}}
\nc{\LC}{\mathcal{I}^{{CINT}}_{_1}\hspace{-0.02in}}
\nc{\QC}{\mathcal{I}^{{CINT}}_{_2}\hspace{-0.02in}}
\nc{\LM}{\mathcal{I}^{{M}}_{_1}\hspace{-0.02in}}
\nc{\QM}{\mathcal{I}^{{M}}_{_2}\hspace{-0.02in}}
\nc{\HMj}{\mathcal{I}^{{M}}_{_{0,j}}\hspace{-0.02in}}

\nc{\wt}{\widetilde}
\renewcommand{\tilde}{\widetilde}
\renewcommand{\bar}{\overline}

\usepackage[titletoc]{appendix}
\usepackage{amssymb,amsfonts,amsmath,amsopn,bm,mathrsfs,graphicx,
color,subfigure,bbm,morefloats}
\usepackage{changes}

\renewcommand{\bf}{\textbf}

\usepackage{fullpage}
\usepackage{bm}
\usepackage{url}
\newtheorem{theorem}{Theorem}

\newtheorem{remark}[theorem]{Remark}
\renewcommand{\theequation}{\thesection.\arabic{equation}}

\begin{document}

\title{Second-Harmonic Imaging in Random Media} 

\author{Liliana Borcea}
\address{Department of Mathematics, University of Michigan,
Ann Arbor, MI 48109}
\email{borcea@umich.edu}

\author{Wei Li}
\address{Department of Mathematics, University of Michigan,
Ann Arbor, MI 48109}
\email{leewei@umich.edu}

\author{Alexander V. Mamonov}
\address{Department of Mathematics, University of Houston, Houston, TX 77004}
\email{mamonov@math.uh.edu}

\author{John C. Schotland}
\address{Department of Mathematics and Department of Physics, University of Michigan,
Ann Arbor, MI 48109}
\email{schotland@umich.edu} 

\begin{abstract}
We consider the problem of optical imaging of small nonlinear scatterers in random media. 
We propose an extension of coherent interferometric imaging (CINT) that applies to scatterers that emit second-harmonic light.  We compare this method  to a nonlinear version of migration imaging and find that the images obtained by CINT are more robust to statistical fluctuations. This finding is supported by a resolution analysis that is carried out in the setting of geometrical optics in random media.  It is also consistent with numerical simulations for which the assumptions of the geometrical optics model do not hold. 
\end{abstract}

\maketitle

\section{Introduction}
\subsection{Background}
There has been considerable recent interest in the development of methods for optical imaging of biological systems in the mesoscopic scattering regime~\cite{arridge2009optical}. There are multiple potential applications incuding imaging of engineered tissues and semitransparent organisms, such as Drosophila and zebra fish, among others~\cite{vinegoni}. Here the term mesoscopic refers to systems whose size is of the order of the transport mean free path of light~\cite{van1999multiple,transport_mean_free_path}. In this setting, light exhibits sufficiently strong scattering so that direct imaging is not possible.
Moreover, imaging modalities that rely on diffuse propagation of light are not effective. Thus, optical imaging techniques that bridge the gap between microscopic and macroscopic scales are of increasing importance. 

Mesoscopic imaging may be carried out using either coherent or incoherent approaches. Incoherent methods, which make use of intensity measurements of transmitted light, include confocal microscopy, optical projection tomography and single-scattering optical tomography~\cite{arridge2009optical,wilson,sharpe,florescu_1,florescu_2,florescu_3,ambartsoumian_1,katsevich,ambartsoumian_2,haltmeier}. Coherent methods, which additionally rely on measurements of the optical phase, include optical coherence microscopy and interferometric synthetic aperture microscopy~\cite{izatt,ralston}. 

An important refinement of optical imaging is to introduce a fluorescent molecular probe that binds to a target of interest. In this manner, spectral isolation and increased signal-to-noise of the detected light can be achieved~\cite{vinegoni}. Spectral isolation may also be realized by utilizing a contrast agent that exhibits a nonlinear optical response. This approach has the advantage that it is not affected by fluorescent photobleaching. Experiments in which second-harmonic generation (SHG) has been utilized for mesoscopic imaging have recently been reported~\cite{brown,hsieh}. The theory of SHG is reviewed in~\cite{Boyd}. Briefly, second-harmonic light is generated by a nonlinear process in which the polarization of a material medium depends both linearly and quadratically on the electric field. 

In this paper we introduce a new form of coherent mesoscopic imaging, in which interferometric measurements  are used to localize small scatterers that emit second-harmonic light. The scatterers are embedded in an unknown background medium, that is taken to be random. The method that is proposed is an extension of coherent interferometry (CINT) to media that exhibit a nonlinear response. The key idea of CINT is to form images by backpropagating correlations of the measured field
\cite{borcea2011enhanced,borcea2006adaptive}. We compare nonlinear CINT to a nonlinear version of migration imaging, in which the fields alone are backpropagated \cite{bleistein2013mathematics}. We find that the images obtained by CINT are more robust to statistical fluctuations in the background medium in comparison to images obtained by migration. This conclusion is supported by an analysis of image resolution carried out within the framework of a model for geometrical optics in random media.  It is also consistent with numerical simulations of wave propagation in media for the which the assumptions of the geometrical optics model do not hold. We note that migration imaging for nonlinear media was also studied in~\cite{Ammari}. However, in that work, only migration imaging was investigated in the setting of a two-dimensional model for SHG with transverse-magnetic polarization for the incident field and transverse-electric polarization for the second harmonic.

\subsection{Physical principles}
\begin{figure}[t]
\centering
\includegraphics[width=0.4\linewidth]{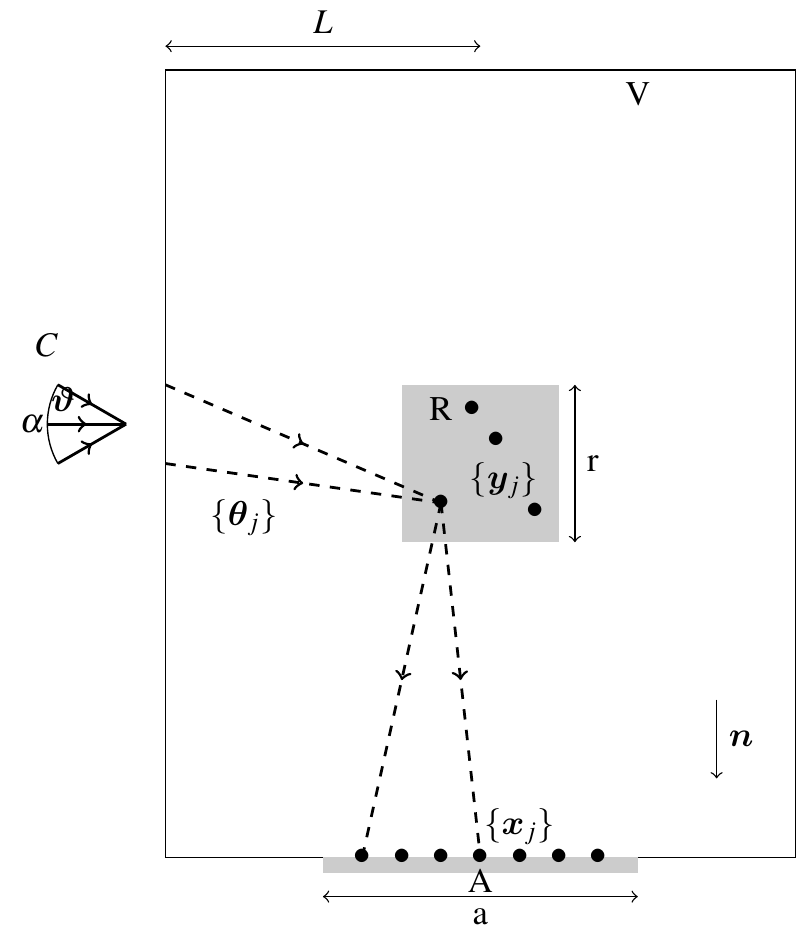}
\caption{Setup of the imaging problem. $A$ denotes the array, $C$ the cone of incident directions, $R$ the search region and $V$ the imaging domain. }
\label{imsetup}
\end{figure}

We are interested in imaging small nonlinear scatterers at positions $\by_j$, for $j = 1, \ldots, N_y$, in a medium occupying a bounded domain 
$V \subset \mathbb{R}^m$, with boundary $\partial V$, for $ m \ge 2$, as   illustrated in
Figure~\ref{imsetup}. We restrict our attention to the case of second harmonic generation (SHG), but more general, quadratic or cubic nonlinearities could be treated similarly.
For convenience, we  let $V$ be a cube (or square in two dimensions) of side $2 L$.
The medium in $V$ is illuminated by monochromatic plane
waves at frequency $\omega$, in the directions of the unit vectors
$\bth_j$, for $j=1,\ldots,N_{\theta}$. These vectors belong to a cone
$C$ with axis along the unit vector $\bvth$ normal to $\partial V$, and small opening angle
$\alpha$.  The resulting waves are measured by an array of detectors located at points $\bx_s$ in the array aperture ${A}$, for
$s=1,\ldots,N_{\bx}$. The array lies on one side of the boundary $\partial
V$,  and $A$ is a square of  side $a \ll L$, or  a segment of length $a \ll L $ in two dimensions.  In the analysis, not the numerics, we suppose that the 
scatterers are confined to a  cubic (square) region $R$ with  the same center as $V$, of side $r \ll L$.  We also let the centers of $A$ and $V$ be lined up  along the unit normal vector
${\bn}$, which  is orthogonal to $\bvth$ and the array aperture.

The imaging problem is to estimate the locations of the  $N_{\by}$ nonlinear
scatterers from the measurements of the wave fields at the array, at frequency $\om$
and the second harmonic $2 \om$. When the
medium in $V$ is known and is non-scattering, we can image with coherent
methods known as matched filtering in the signal processing
\cite{therrien1992discrete} and radar literature
\cite{nolan2002synthetic}, migration in seismic imaging
\cite{biondi20063d}, and backprojection in computed tomography
\cite{natterer1986mathematics}. These methods form an image by superposing the array
measurements backpropagated in the known medium to imaging points $\bys$
in $R$. The backpropagation is done analytically when
the Green's function is known, or numerically, and its purpose is to
compensate for the phases of the measurements at points $\bys$ near the
scatterer location, so that they add constructively, and the
imaging function displays a peak. We refer henceforth to such imaging
as migration.

In this paper we consider the problem of imaging in heterogeneous media, containing
numerous small inhomogeneities. A single
inhomogeneity by itself is assumed to be much weaker than the scatterer that
we wish to image. However, the waves interact with many inhomogeneities over the $O(L)$ distances of propagation, 
and this results in strong cumulative scattering effects that must be taken into account in imaging. 
This problem is challenging because, in applications, it is impossible to know the numerous  inhomogeneities, and we cannot hope
to determine them from the data. Thus, imaging is carried in an
uncertain environment. We incorporate this uncertainty in the data
model by studying imaging in random media. The goal is to understand
whether it is possible to obtain robust estimates of the nonlinear
scatterer locations from the measurements registered by the array in a fixed
realization of the random medium. Such robustness is known as
statistical stability, and it means that the images vary little from one realization of the random medium to another.

Our results consist of analysis and numerical simulations of imaging of nonlinear
scatterers  in a random medium. We consider two regimes where cumulative
scattering by the inhomogeneities causes large distortion of the wave
field measured at the array. In physical terms, this means that $L$ is larger than the scattering length (also known as scattering mean free path) ${\it \ell^s}$,
which is the characteristic length scale on which the waves randomize
\cite{van1999multiple}. The analysis uses a geometrical optics model
\cite{rytov1989principle}, where the typical size $\ell$ of the
inhomogeneities is large with respect to the wavelength $\la$, and $L
\gg {\it \ell^s} \gg \ell$, so that the net scattering effects in the medium
amount to large random wavefront distortions.  In the numerics we
consider a regime with $\ell \sim \la$, where the waves interact more
efficiently with the inhomogeneities. Because we are interested in
coherent imaging, we take $L$ of the order of the transport mean free
path, which is the characteristic distance at which the waves forget
their initial direction due to scattering~\cite{van1999multiple}. As explained above, this corresponds to the mesoscopic regime. At larger distances, the waves are in a radiative transfer regime, and only incoherent imaging methods are effective~\cite{arridge2009optical}.

We find in both the analysis and numerical simulations that the statistical
fluctuations of the array measurements due to scattering in the random medium are very different than those of additive
noise. They cannot be mitigated simply by summation of the measured fields, as in
migration. Ideally, the fluctuations could be treated by backpropagation of the measurements in the  medium, as in time
reversal \cite{fink1992time}. However, this cannot be done
because the medium is unknown.  We can only backpropagate in a
hypothetical reference medium with known wave velocity, which is constant in this paper, and thus obtain  migration 
images. We show that these images are statistically unstable, meaning that they change significantly and unpredictably from one realization of the random medium to another.

The coherent interferometric (CINT) method
\cite{borcea2006adaptive,borcea2007asymptotics} is designed to
ameliorate the effects of random wave distortions, by imaging with cross-correlations of the measurements instead of the measurements themselves.  The CINT
imaging function is given by the superposition of cross-correlations
backpropagated in the reference medium. Its mathematical expression resembles that of the time reversal function analyzed in
\cite{papanicolaou2004statistical,papanicolaou2007self}, and its
statistical stability and resolution are studied in
\cite{borcea2007asymptotics,borcea2011enhanced}. Unlike in time
reversal, where super-resolution of focusing occurs, the stability of
CINT comes at the expense of resolution, which is determined by two
characteristic scales in the random medium: the decoherence frequency
and length. These scales describe the frequency offsets and detector
separations over which the waves become statistically uncorrelated,
and must be taken into account in the calculation of the
cross-correlations \cite{borcea2006adaptive}. Moreover, they determine
the conditions under which CINT is statistically stable. This can be
formally understood as a consequence of the law of large numbers, due
to the summation in the imaging function of many statistically
uncorrelated terms, when the array aperture  is much larger than
the decoherence length, and the probing signals have bandwidth that is
larger than the decoherence frequency \cite{borcea2007asymptotics}.

In this paper we study CINT imaging with time harmonic waves, so summation over frequencies is not carried out. Such a summation is essential for
the statistical stability of CINT, and of time reversal for that
matter \cite{blomgren2002super}, when scattering in the medium causes
significant delay spreading.  Here we consider weaker
scattering regimes, like  random geometrical optics, where imaging can
be performed at a single frequency \cite{borcea2011enhanced}.

The remainder of this paper is organized as follows. In Section 2, we present the formulation of the problem and introduce the necessary migration and CINT imaging formalism. In Section 3, we analyze the migration and CINT point spread functions in the setting of the random geometrical optics model. Section 4 presents our numerical results for two-dimensional model systems. We end with a summary in Section 5. The appendices present some technical details and derivations of results 
presented in the paper.

\section{Formulation of the problem}
\label{sect:form}

In this section we describe the model of the imaging experiment and define the
migration and CINT imaging functions that are analyzed theoretically and
numerically in the following sections.

\subsection{Random model of the array data}
\label{sect:form.1}
We consider for simplicity scalar waves, that obey the Helmholtz
equation with a random wave speed $c(\bx)$ of the form
\begin{equation}
c(\bx) = c_0 \left[ 1 + 4 \pi \eta(\bx)\right]^{-1/2},
\label{eq:f1}
\end{equation}
where $c_0$ is constant and $\eta(\bx)$
is the linear susceptibility of the medium. Here $\eta(\bx)$ is taken to be a mean zero, stationary random process that is bounded almost surely, so that the right hand
side in \eqref{eq:f1} remains positive. We also assume that
$\eta(\bx)$ is mixing (lacks long range correlations) \cite{kushner}, with integrable autocorrelation function. 
The small
scatterers embedded in the random medium with wave speed \eqref{eq:f1} are modeled by  the linear
susceptibility $\eta_{1}(\bx)$ and quadratic susceptibility
$\eta_{2}(\bx)$. These are functions with small amplitude and support 
concentrated in the vicinity of the points $\by_j$, for $j = 1, \ldots N_{\by}$.

Under conditions of weak nonlinearity, the linear and second harmonic wave fields obey~\cite{Boyd}
\begin{align}
\label{eq:waveu}
\Delta u_1(\bx;\bth)+k^2 [1+4\pi\eta(\bx)]u_1( \bx;\bth)&=-4\pi
k^2[\eta_1(\bx)u_1( \bx;\bth)+2\eta_2(\bx)u_2(
  \bx;\bth)u_1^*(\bx;\bth)]\ ,\\ \Delta u_2( \bx;\bth)+(2k)^2
[1+4\pi\eta(\bx)]u_2( \bx;\bth)&=-4\pi (2k)^2 [ \eta_1(\bx)u_2(
  \bx;\bth)+\eta_2(\bx)u_1^2( \bx;\bth)],
\label{eq:wavev}
\end{align}
for $\bx \in V$. Here $u_1(\bx;\bth)$ is the wave
at frequency $\om$ and $u_2(\bx;\bth)$ is the wave at the second
harmonic $2 \om$. 
The medium is illuminated by $N_{\theta}$ plane waves at frequency $\om$ 
\begin{equation}
u_{1,0}^{(i)}(\bx;\bth_q) = e^{i k \bth_q \cdot \bx}, 
\label{eq:uinc}
\end{equation}
where $k = \om/c_0$ is the wavenumber and $\bth_q$ are unit wave vectors. The scattered waves, given by
$u_1(\bx;\bth_q)-u_{1,0}^{(i)}(\bx;\bth_q)$ and $u_2(\bx;\bth_q)$,
satisfy outgoing boundary conditions on $\partial V$.

The imaging is to be carried using measurements of the wave fields at the $N_{\bx}$ detectors in the array, gathered in the set
\begin{equation}
\mathcal{U} = \left\{ u_1(\bx_s;\bth_q),~ u_2(\bx_s;\bth_q), ~ ~~s = 1,
\ldots, N_{\bx}, ~ ~ q = 1, \ldots, N_{\theta}\right\}.
\label{eq:datarand}
\end{equation}
Note that $u_1(\bx;\bth_q)$ and $u_2(\bx;\bth_q)$ are random fields, and the actual
array measurements are for a single realization of the medium. These
are the data used to form images, as explained in the next
sections. We only use the set $\mathcal{U}$ for the statistical
analysis of the imaging functions.

\subsection{Migration imaging}
\label{sect:form.2}
In migration imaging  one assumes that the background medium is known and non-scattering, and that the waves are reflected  once at the unknown scatterers. 
In our setting this corresponds to assuming that  $\eta \equiv 0$ in \eqref{eq:f1},  and 
using the Born approximation with respect to $\eta_1(\bx)$ and $\eta_2(\bx)$ in  equations \eqref{eq:wavev}. We obtain  the linearized forward mappings from the unknown susceptibilities 
$\eta_1(\bx)$ and $\eta_2(\bx)$  to the scattered fields at the array
\begin{equation}
 F_{j}[\eta_{j}](\bx_s,\bth_q)= (jk)^2\int_V d \by \,\eta_{j}(\by)
 G_0(\bx_s,\by;j\omega)e^{ijk\bth_q\cdot\by}, \qquad j = 1, 2,
\label{eq:maphom}
\end{equation}
where $G_0$ is the outgoing Green's function, given by
\begin{equation}
\label{Green03}
G_0(\bx,\by;\omega)=\frac{e^{ik|\bx-\by|}}{|\bx-\by|} 
\end{equation}
in three dimensions, and by
\begin{equation}
\label{Green02}
G_0(\bx,\by;\omega)=i\pi H_0^{(1)}(k|\bx-\by|)
\end{equation}
in two dimensions, where $H_0^{(1)}$ is the Hankel function of the
first kind. The mappings in \eqref{eq:maphom} are evaluated at  the detector locations   indexed by $s = 1, \ldots, N_{\bx}$, 
and at the incident plane wave vectors 
indexed by $q = 1, \ldots,N_{\theta}$. The index $j = 1, 2$ corresponds to the frequencies by  $j  \om$.  

Let us denote by $d_1$ and $d_2$ the data for inversion,  given by
\begin{align}
d_1(\bx_s,\bth_q) &\approx
u_{_{1}}^{(real)}(\bx_s;\bth_q)-u_{1,0}^{(i)}(\bx_s;\bth_q) , \label{eq:data1rand}\\ d_{_{2}}(\bx_s,\bth_q) 
&\approx u_{_{2}}^{(real)}(\bx_s;\bth_q) , \quad
s = 1, \ldots, N_{\bx}, ~~ q = 1,\ldots, N_{\theta},
\label{eq:data2rand}
\end{align}
where $u_j^{(real)}$, for $j = 1,2 $, denote the measurements in the real medium, one realization of the random model.
The approximation sign accounts for noise that is unavoidable in experiments. When this noise is additive, 
mean zero, identically distributed, with finite variance, the Gauss-Markov theorem gives that the best unbiased 
linear estimator of the unknown susceptibilities is the solution of the least squares problem  \cite{Rao}
\begin{equation}
\eta_{j}^{^{LS}} = \underset{\eta_{j}}{\mbox{arg min}}\, \|d_{j} -
F_{j}[\eta_{j}] \|_2^2 \label{eq:LS}.
\end{equation}
The minimizers  $\eta_{j}^{^{LS}}$  satisfy the normal equations 
\begin{equation}
F_{j}^H F_{j}[\eta_{j}^{^{LS}}] (\bys) =
F_{j}^H [d_{j}](\bys) = (jk)^2 \sum_{s=1}^{N_{\bx}}
\sum_{q=1}^{N_{\theta}} G_0^\star(\bys,\bx_s;j\om) e^{-i j k
  \bth_q \cdot \bys} d_{j}(\bx_s,\bth_q),
\label{eq:LS1}
\end{equation}
where the superscript $H$ denotes the adjoint with respect to the Euclidian inner
product. The integral (normal)
operators $F_{j}^H F_{j}$ map
$\eta_{j}^{^{LS}}$ to
\begin{equation}
F_{j}^H F_{j}[\eta_{j}^{^{LS}}] (\bys) = \int_V d \by \, \eta_{j}^{^{LS}}(\by) \mathcal{K}_{j}(\bys,\by), 
\label{eq:LS2}
\end{equation}
and their integral kernels
\begin{equation}
\mathcal{K}_{j}(\bys,\by) = (jk)^4\sum_{s=1}^{N_{\bx}} \sum_{q=1}^{N_{\theta}}
G_0(\bx_s,\by;j\om) G^\star_0(\bys,\bx_s;j\om) e^{i j k \bth_q
  \cdot (\by-\bys)} 
\label{eq:LS3}
\end{equation}
are equal, up to constant factors, to the time reversal point spread functions
at frequencies $j \om$, for $j = 1, 2$. 
In our setting, the kernels \eqref{eq:LS3} peak along the diagonal
$\bys = \by$, so we can approximate the left hand side in
\eqref{eq:LS2} by $\eta_{j}^{^{LS}}(\bys)$ multiplied by a
constant. 

In imaging, we are interested in the support of the
sussceptibilities, so we can neglect the constants and obtain from
\eqref{eq:LS3} the migration imaging functions
\begin{equation}
\IMj(\bys) = \sum_{s=1}^{N_{\bx}} \sum_{q=1}^{N_{\theta}}
G_0^\star(\bys,\bx_s;j\om) e^{-i j k \bth_q \cdot \bys}
d_{j}(\bx_s,\bth_q), \quad j = 1,2.
\label{eq:LS4}
\end{equation}
Since $k |\bys - \bx_s| \gg 1$ and the array aperture size $a$ is much
smaller than $|\bx_s-\bys| = O(L)$, for $\bys$ in the search region $R$
(recall Figure~\ref{imsetup}), we can approximate the Green's function
in \eqref{eq:LS4} by
\begin{equation}
G_0(\bys,\bx_s;j\om) \approx C_j e^{i j k |\bys-\bx_s|}, \label{eq:LS5}
\end{equation}
for constant $C_j$, with $j = 1, 2$. This result holds in both two and three dimensions, as follows from the
asymptotics of the Hankel function in \eqref{Green02}. Thus, the right
hand side of \eqref{eq:LS4} is the superposition of the measurements,
with phases given relative to the imaging point $\bys$.  The
superposition is needed for focusing the image and averaging over the
noise. When $\bys$ is close to a scatterer location, and the medium is either homogeneous or has negligible effect on the waves, so that their propagation is approximated by $G_0$, the phases in $d_{j}$ are cancelled approximately.  Then, the  terms in \eqref{eq:LS4} add constructively and the imaging function displays a peak. 
In this paper we are interested in imaging in stronger scattering media, where $G_0$ is not a good 
model for wave propagation,  and  the migration imaging function \eqref{eq:LS4} either does not focus,
or gives spurious peaks at locations that may not be close to the scatterers.

\subsection{Coherent interferometric imaging}
\label{sect:form.3}

Let $\bys$ be an imaging point and define 
\begin{equation}
b_{j}(\bx_s,\bth_q,\bys) = d_{j}(\bx_s,\bth_q)
G_0^\star(\bys,\bx_s;j\om) e^{-i j k \bth_q \cdot \bys},\label{eq:CINT1}
\end{equation}
for detector index $s =
1, \ldots,N_{\bx}$, plane wave index $q = 1, \ldots, N_{\theta}$, and frequency index $j = 1, 2.$
The CINT imaging function is formed by superposition of local
cross-correlations of $b_{j}$. By local we mean that 
we cross-correlate only at nearby detectors and for nearby incident directions
\begin{equation}
|\bx_s - \bx_{s'}| \le X, \quad |\bth_q - \bth_{q'}| \le \Theta,
\label{eq:CINT2}
\end{equation}
where $\Theta$ and $X$ are scales that account for the decorrelation
of the waves due to scattering in the random medium
\cite{borcea2006adaptive}. Intuitively, in the language of geometrical optics, waves traveling along different paths are decorrelated because they interact with different parts of the random
medium, assumed to lack long range correlations of the fluctuations of the wave speed.  Note that
in practice the scales of decorrelation of the waves  are usually unknown, so they must be estimated from the data, either
using statistical data analysis or by optimization, which seeks to
improve the focusing of CINT images, as explained in
\cite{borcea2006adaptive}.  We denote henceforth the true decorrelation parameters in
the medium by $X_{d,j}$ and $\Theta_d$, to distinguish them from those
used in the calculation of the cross-correlations, and assume that
\begin{equation}
{X}/{X_{d,j}} = O( 1), \quad {\Theta}/{\Theta_{d}} = O(1).
\label{eq:decorestim}
\end{equation}
We also refer to $X_{d,j}$ as decoherence lengths and $\Theta_d$ as decoherence angle. Note that 
the decoherence length is proportional to the wavelength, so it depends on $j$. We suppress, for simplicity of notation, the 
dependence of the thresholding parameter $X$ on $j$.

Let us introduce the center and offset sensor locations
\begin{equation}
\obx_{ss'} = ({\bx_s + \bx_{s'}})/{2}, \quad \tbx_{ss'} = \bx_s - \bx_{s'},
\quad s,s' = 1, \ldots, N_{\bx},
\label{eq:CINT3}
\end{equation}
and direction vectors
\begin{equation}
\oth_{qq'} = ({\bth_q + \bth_q'})/{2}, \quad \tth_{qq'} = \bth_q -
\bth_{q'}, \quad q,q' = 1, \ldots, N_{\theta}.
\label{eq:CINT4}
\end{equation}
We count the center variables by $\bar{s}= 1, \ldots,
\bar{N}_{_{\bx}}$ and $\bar{q} = 1, \ldots, \bar{N}_{_\theta}$, and
the offsets by $\tilde{s} = 1, \ldots, \tilde{N}_{_{\bx}}$ and
$\tilde{q} = 1, \ldots, \tilde{N}_{_\theta}$.  The local
cross-correlations are
\begin{equation}
\mathcal{C}_j(\obx_{\bar{s}},\oth_{\bar{q}},\bys) = \sum_{\tilde{s} =
  1}^{\tilde{N}_{_{\bx}}} \sum_{\tilde{q}=1}^{\tilde{N}_\theta}
\Phi \Big( \frac{\tbx_{\tilde{s}}}{X}\Big)
\Phi \Big(\frac{\tth_{\tilde{q}}}{\Theta}\Big)
b_j \Big(\obx_{\bar{s}} + \frac{\tbx_{\tilde{s}}}{2},
\oth_{\bar{q}} + \frac{\tth_{\tilde{q}}}{2}, \bys\Big)
b_j^\star  \Big(\obx_{\bar{r}} -
  \frac{\tbx_{\tilde{s}}}{2}, \oth_{\bar{q}} -
  \frac{\tth_{\tilde{q}}}{2}, \bys\Big),
\label{eq:CINT5}
\end{equation}
where $\Phi$ is a smooth window  of support of order one, used to limit the 
detector and director offsets  by $X$ and  $\Theta$.
The CINT imaging function is formed by the superposition of
\eqref{eq:CINT5}, 
\begin{equation}
\ICINTj(\bys) = \sum_{\bar{s}=1}^{\bar{N}_s}
\sum_{\bar{q}=1}^{\bar{N}_\theta}
\mathcal{C}_j(\obx_{\bar{s}},\oth_{\bar{q}},\bys).
\label{eq:CINT6}
\end{equation}
Were it not for the windows $\Phi$ in \eqref{eq:CINT5}, this expression would equal the 
square of the migration imaging function \eqref{eq:LS4}. The windows play a smoothing role 
in \eqref{eq:CINT6}, by  convolution. This blurs the images, but is essential for stabilizing them statistically  \cite{borcea2007asymptotics}.

\begin{remark}
{\rm 
The phase compensation in equation ~\eqref{eq:CINT1} assumes the complete removal of the direct 
waves that have not interacted with the scatterers that we wish to
image.  In homogeneous media this is achieved by the subtraction of the incident wave $u_{1,0}^{(i)}$ from the measurements.
However, this 
wave is affected by scattering in random media,  so the
subtraction in \eqref{eq:data1rand} does not achieve its purpose. The
unwanted direct wave may be removed in our geometrical setting if the array of detectors can differentiate among arrival directions,
 and the scattering medium is  not strong enough to mix the directions of the waves, as is the case in
the geometrical optics regime considered in our analysis. Such
differentiation may be achieved for example by an approximate plane
wave decomposition of the measurements, using a discrete Fourier
transform with respect to the coordinates of the detectors in the  array.  We do not perform such a differentiation here and work instead
with \eqref{eq:data1rand}, to illustrate the  effect of the
unwanted direct arrivals on the imaging of $\eta_1$. This problem
does not extend to the second harmonic wave \eqref{eq:data2rand}, which is emitted at the
nonlinear scatterers.}
\end{remark}

\section{Analysis of the migration and CINT point spread functions}
\label{sect:anal}
\setcounter{equation}{0}
In this section we analyze migration and CINT imaging of a single scatterer at location $\by$ in $V$. This means that we analyze the point spread functions of \eqref{eq:LS4} and \eqref{eq:CINT6}. We consider a geometrical optics wave propagation regime, with large random distortions of the wavefront. The analysis is in most respects the same in two and three dimensions, so we focus our attention on the three dimensional case. 

To begin, we write the susceptibility of the medium as
\begin{equation}
4 \pi \eta(\bx) = \sigma \mu\left(\frac{\bx}{\ell}\right),
\label{eq:A1}
\end{equation}
where $\mu$ is a random, stationary process with mean zero and Gaussian autocorrelation
\begin{equation}
\EE[\mu(\bm{h})\mu(\bm{0})] = e^{-\frac{|\bm{h}|^2}{2}}.
\label{eq:A2}
\end{equation}
This choice of autocorrelation is not neccesary, but it convenient for the analysis, because it allows us to obtain explicit expressions for the statistical moments of the imaging
functions.  The process $\mu$ is normalized so that the maximum of  \eqref{eq:A2} 
equals one, and
\[
\int_{\mathbb{R}^3} d \bx \,
\EE\left[\mu\left(\frac{\bx}{\ell}\right)\mu(\bm{0})\right] = (2
\pi)^{3/2} \ell^3.
\]
Thus, the scale $\sigma$ in \eqref{eq:A1} characterizes the amplitude of
the random fluctuations of the susceptibility, and $\ell$ characterizes the
correlation length.

In section \ref{sect:anal.0} we introduce the necessary scalings and then
describe in section \ref{sect:anal.1} the random geometrical optics
model of wave propagation. We base our analysis on the linearized data model
defined in section \ref{sect:LINmod}, justified by the weak
nonlinearity. With this model, we calculate in sections
\ref{sect:anal.2} and \ref{sect:anal.3} the expectation and variance
of the migration and CINT images, in order to study their resolution
and statistical stability. 

\subsection{Scaling}
\label{sect:anal.0}
We assume for convenience that the scatterer location  $\by$ is at  the center of $V$, and that the center of the array aperture $A$ is lined up with  $\by$, along the unit normal $\bn$ to the boundary that supports the array. We also suppose that $\bn$ is orthogonal to the unit vector $\bvth$ along the axis of the cone of illuminations. Since the aperture size $a$ is much smaller than $L$, and the cone of incident 
directions $\bth_q$ has a small opening angle, we have 
\[
|\bx_s-\by| \approx L, \quad |\by-\by^{(i)}(\bth_q)| \approx L, \]
for $s = 1, \ldots, N_{\bx}$ and  $q = 1, \ldots,
N_{\theta}.$
Here we denote by  $\by^{(i)}(\bth_q)$  the
incident point on $\partial V$ of the ray entering the domain in the
direction $\bth_q$ and passing through $\by$.  
We take the origin of the system of coordinates 
at the center of $A$, with one axis parallel to $\bn$, so that we can write henceforth
$
\bx = (\bx_\perp,0), 
$
for the points in the array, 
with two dimensional vectors $\bx_\perp$ in the aperture $A$. The set of these points is denoted by
\begin{equation}
\label{eq:ap3D}
\mathbb{A} = \left\{ (\bx_\perp,0), \quad \bx_\perp \in A \right\}.
\end{equation}

The random geometrical optics wave propagation model described in the
next section applies to the regime of separation of scales
\begin{equation}
\la \ll \ell \ll L,
\label{eq:A3}
\end{equation}
with small amplitude $\sigma$ of the fluctuations of the
susceptibility, satisfying
\begin{equation}
\sigma \ll 
\left({\ell}/{L}\right)^{3/2}, \quad \sigma \ll {\sqrt{\la \ell}}/{L}.
\label{eq:A4}
\end{equation}
As shown in \cite[Chapter 1]{rytov1989principle}, the first bound in
\eqref{eq:A4} is needed so that the waves propagate along straight
rays, and the variance of the amplitude of the Green's function is
negligible. The second bound ensures that the second order (in
$\sigma$) corrections of the travel time are negligible.  We estimate
in the next section that the standard deviation of the random travel
time fluctuations is of order $\sigma \sqrt{\ell L}/c_0$. When this is small, simpler methods like 
migration work well. We are interested in the case of  
travel time fluctuations that are larger than the  period $2 \pi/\om$. This occurs when 
\begin{equation}
 \sigma \gg \frac{\la}{\sqrt{\ell L}},
\label{eq:A12}
\end{equation}
and  is consistent with \eqref{eq:A4} when $ \ell \gg \sqrt{\la L}$.

We already stated that $a \ll L$. To simplify the calculations, we consider a paraxial
regime\footnote{Imaging may be done with larger apertures and wider
  opening angles of the illumination cone $C$, but the
  expressions of the imaging functions become complicated. The
  analysis presented here may be used in such cases, after segmenting
  the aperture and illumination cone in subsets satisfying our
  assumptions. The results apply for each subset, and the images are
  obtained by summation over the subsets.}  with
\begin{equation}
a \ll (\la L^3)^{1/4} \ll L.
\label{eq:APA1}
\end{equation}
The illumination directions  belong to the cone $C$ with axis along the unit vector
$\vartheta$ assumed orthogonal to $\bn$, and with 
opening angle  $\alpha$   satisfying
\begin{equation}
\alpha = O\left(\frac{a}{L}\right).
\label{eq:APA2}
\end{equation}

The search region $R$ is centered at $\by$. It is a cube of side length $r$ satisfying 
\begin{equation}
r \ll \frac{ \la L^2}{a^2} \ll a,
\label{eq:APAR}
\end{equation}
so that we can  use the following approximation of the Green's function in the reference medium
\begin{equation}
G_0(\bx,\bys;j \om) \approx \frac{1}{L} \exp \left[ i j k
  \left(\ys_{\parallel} + \frac{|\bx_\perp - \bys_{\perp}|^2}{2 L} \right)
  \right], \quad  \, \bys \in R.
\label{eq:APAR3}
\end{equation}
Here we wrote  $\bys = (\bys_\perp,\ys_\parallel)$,  with $\ys_\parallel$ equal to  the distance of $\bys$ from the array, along
$\bn$, and $\bys_\perp$ the two dimensional vector in the plane orthogonal to $\bn$. With this notation we note that $y_\parallel = L$.

We expect from the analysis in \cite{borcea2011enhanced} that to
obtain statistically stable CINT images we need $a \gg \ell$. There
are many scalings that allow $\sqrt{\la L} \ll \ell \ll a$, so we
choose one that simplifies slightly the moment calculations of the
random travel time corrections. Specifically, we consider the length
scale ordering
\begin{equation}
\la \ll \sqrt{\la L} \ll \ell \ll (\la L^2)^{1/3} \ll a \ll (\la L^3)^{1/4} \ll L,
\label{eq:as1}
\end{equation}
and gather the assumptions \eqref{eq:A4}-\eqref{eq:A12} on $\sigma$ in 
\begin{equation}
\frac{\la}{\sqrt{\ell L}} \ll \sigma \ll \frac{\sqrt{\la \ell}}{L}.
\label{eq:as2}
\end{equation}
Here we used that\footnote{Note  that \eqref{eq:as1} is consistent because
$\displaystyle
\frac{\sqrt{\la L}}{(\la L^2)^{1/3}} =
  \left(\frac{\la}{L}\right)^{1/6} \ll 1$ and $\displaystyle \frac{(\la
    L^2)^{1/3}}{(\la L^3)^{1/4}} = \left(\frac{\la}{L}\right)^{1/12}
    \ll 1.$}

\[
\frac{(\ell/L)^{3/2}}{\sqrt{\la \ell}/L} = \frac{\ell}{\sqrt{\la L}} \gg 1.
\]

\subsection{Random geometrical optics model}
\label{sect:anal.1}
We refer to \cite[Chapter 1]{rytov1989principle} and \cite[Appendix
  A]{borcea2011enhanced} for the derivation of the geometrical optics
model. It holds in the scaling regime defined by equations
\eqref{eq:as1}-\eqref{eq:as2}.

The geometrical optics approximation of the Green's function, denoted
by $G$, is
\begin{equation}
G(\bx,\by; j \om) = G_0(\bx,\by;j \om) e^{i j k \nu(\bx,\by)}, \quad  \, \bx \in \mathbb{A}, 
\label{eq:A5}
\end{equation}
where $G_0$ is given by \eqref{Green03}, and the random phase $\nu$ is given by the 
integral of the fluctuations $\mu$ along straight rays
\begin{equation}
\label{eq:A6}
\nu(\bx,\by)=\frac{\sigma|\bx-\by| }{2}\int_0^1dt \, \mu
\Big(\frac{(1-t)\by}{\ell}+\frac{t\bx}{\ell}\Big).
\end{equation}
The approximation of the direct wave, which enters the medium as the
plane wave \eqref{eq:uinc}, is
\begin{equation}
u_1^{(i)}(\bx;\bth) = e^{i k \bth \cdot \bx + i k \gamma(\bx,\bth)},
\label{eq:A7}
\end{equation}
with random phase 
\begin{equation}
\gamma(\bx,\bth)=\frac{\sigma|\bx - \bx^{(i)}(\bth)|}{2}\int_0^1dt\,
      {\mu}\Big(\frac{(1-t)\bx}{\ell}+\frac{t\bx^{(i)}(\bth)}{\ell}\Big).
\label{eq:A8}
\end{equation}
Because  
\[|\bx-\by| \approx L , \quad   \bx \in \mathbb{A},
\]
and 
\[
|\bx-\bx^{(i)}(\bth)| = O(L), \quad  \bx \in \mathbb{A} \cup \{\by\},
\]
we can use \cite[Lemma 3.1]{borcea2011enhanced} to conclude that the normalized processes
\begin{equation}
\tilde{\nu}(\bx,\by) = \frac{2}{(2 \pi)^{1/4}}
\frac{\nu(\bx,\by)}{\sigma \sqrt{\ell |\bx-\by|}}, \quad
\tilde{\gamma}(\bx,\bth) = \frac{2}{(2 \pi)^{1/4}}
\frac{\gamma(\bx,\bth)}{\sigma \sqrt{\ell |\bx-\bx^{(i)}(\bth)|}},
\label{eq:A9}
\end{equation}
converge in distribution to Gaussian ones in the limit $\ell/L \to 0$.  

The processes in
\eqref{eq:A9} are mean zero, with variance
\begin{align}
\EE[\tilde{\nu}^2(\bx,\by)] &= \frac{|\bx-\by|}{\ell \sqrt{2 \pi}}
\int_0^1 dt \int_0^1 dt' \exp\left[ -\frac{(t-t')^2|\bx-\by|^2}{2 \ell^2}\right]
\approx 1, \\ \EE[\tilde{\gamma}^2(\bx,\bth)] &=
\frac{|\bx-\bx^{(i)}(\bth)|}{\ell \sqrt{2 \pi}} \int_0^1 dt \int_0^1 dt'
\exp\left[ -\frac{(t-t')^2|\bx-\bx^{(i)}(\bth)|^2}{2 \ell^2}\right] \approx 1,
\end{align}
where we have used \eqref{eq:A2}.
Thus, the variance of the random phases in
\eqref{eq:A5} and \eqref{eq:A7} is
\begin{align}
k^2\EE[{\nu}^2(\bx,\by)] &\approx \frac{\sqrt{2 \pi} \sigma^2 k^2 \ell
  |\bx-\by|}{4} = O\left(\sigma^2 \frac{\ell L}{\la^2}
\right), \label{eq:A10} \\ k^2\EE[{\gamma}^2(\bx,\bth)] &\approx
\frac{\sqrt{2 \pi} \sigma^2 k^2 \ell |\bx-\bx^{(i)}(\bth)|}{4} =
O\left(\sigma^2 \frac{\ell L}{\la^2} \right), \label{eq:A11}
\end{align}
and we conclude from the assumption \eqref{eq:as2} that cumulative
scattering in the medium has a significant net effect on the waves,
manifested as large random wavefront distortions.
\subsubsection{Randomization of the waves}
\label{sect:anal.1.1}
Because the processes \eqref{eq:A9} are approximately Gaussian in our scaling, we can approximate the expectation of the Green's
function \eqref{eq:A5} by
\begin{align}
\EE \left[ G(\bx,\by;j\om)\right] &= G_0(\bx,\by;j\om)
\EE\left[ \exp\left[i j k \nu(\bx,\by) \right]\right] \nonumber
\\ &\approx G_0(\bx,\by;j\om) \exp \left[ - \frac{(j k)^2
    \EE[\nu^2(\bx,\by)]}{2} \right] \nonumber \\&=
G_0(\bx,\by;j \om)
\exp\left[-\frac{|\bx-\by|}{\ell_j^s}\right], \label{eq:A14}
\end{align}
where $\ell_j^s$ are the scattering lengths
\begin{equation}
\ell_j^s = \frac{8}{\sqrt{2 \pi} \sigma^2 (jk)^2 \ell}, \quad j = 1, 2.
\label{eq:A15}
\end{equation}
The scaling relation \eqref{eq:A12} ensures that $\ell_j^s \ll
|\bx-\by| \approx L$, so the mean Green's function is
exponentially small.  Clearly,
\[
\left| G(\bx,\by;j \om) \right| = \left| G_0(\bx,\by;j \om) \right|,
\]
so the standard deviation  is 
\begin{equation}
\mbox{std}\left[G(\bx,\by;j \om)\right] = \sqrt{ \left| G_0(\bx,\by;j
  \om) \right|^2 - \left|\EE\left[ G(\bx,\by;j\om)\right]\right|^2}
  \approx \left| G_0(\bx,\by;j \om) \right|.
\end{equation}
This is much larger than \eqref{eq:A14}, so the wave is randomized 
by scattering in the random medium. In our regime the randomization  arises due to the large phase $j k \nu$ in 
\eqref{eq:A5}.

A similar calculation for the direct wave \eqref{eq:A7} gives 
\begin{align}
\EE \left[u_1^{(i)}(\bx;\bth)\right] &= e^{i k \bth \cdot \bx}
\EE\left[ \exp\left[i k \gamma(\bx,\bth) \right]\right] \nonumber
\\ &= \exp \left[ {i k \bth \cdot \bx}- \frac{k^2
    \EE[\gamma^2(\bx,\bth)]}{2} \right] \nonumber \\&=
\exp\left[i k \bth \cdot \bx -\frac{|\bx-\bx^{(i)}(\bth)|}{\ell_1^s}\right], \label{eq:A16}
\end{align}
and since $|\bx-\bx^{(i)}(\bth)| = O(L) \gg \ell^s_1$ and
$
\big|u_1^{(i)}(\bx;\bth)\big| = 1,
$
we conclude that $u_1^{(i)}(\bx;\bth)$ is randomized and therefore very
different than the incident plane wave \eqref{eq:uinc}, for $\bx \in \mathbb{A}$ and $\bth \in C$.

\subsubsection{Decorrelation of the waves}
\label{sect:anal.1.2}
The statistical  moments of the
wave fields are determined by the second moments of the phases \eqref{eq:A6} and \eqref{eq:A8}, which are 
approximately Gaussian. These
moments are derived in appendix \ref{app:secondmom},  using the
assumptions \eqref{eq:as1}-\eqref{eq:as2}. We 
use them in appendix \ref{ap:lem2} to derive the second moments of the wave fields stated here.

First, let $\bx,\bx'$ be two points in  $\mathbb{A}$. The second moments of the
Green's function \eqref{eq:A5} are 
\begin{align}
\EE\left[ G(\bx,\by;j\om){G^\star(\bx',\by;j
    \om)}\right] &\approx \frac{1}{L^2} \exp\left[i j k \left(
  L + \frac{|\bx_\perp-\by_{\perp}|^2}{2 L}\right)-
  \frac{|\bx'_\perp-\bx_\perp|^2 }{2 X_{d,j}^2}\right],
\label{eq:MOMG}
\end{align}
where we recall that $\bx_\perp$ and $\by_{\perp}$ are the
components of $\bx$ and $\by$ in the plane orthogonal to $\bn$,
and $L$ is the distance from $\by$ to the center of the array.
The length scales $X_{d,j}$ are given by
\begin{equation}
X_{d,j} = \ell \sqrt{\frac{3 \ell_j^s}{2 L}} = O\left(\frac{\la
  \sqrt{\ell}}{\sigma \sqrt{L}}\right) \ll
\ell.
\label{eq:defXd}
\end{equation}
The decay in \eqref{eq:MOMG} with the detector offsets models the decorrelation 
of the waves due to scattering 
in the random medium, so we call $X_{d,j}$ the decoherence lengths. By 
equation \eqref{eq:as1}, we have $X_{d,j} \ll a$, which is essential for obtaining statistically 
stable CINT images, as we show later.

Next, let $\bx,\bx'$ be two points in $\mathbb{A}$, and $\bth$ and $\bth'$ 
be two illumination directions in the cone $C$. The second moments of the direct wave are 
\begin{align}
\EE\left[u_1^{(i)}(\bx,\bth) {u_1^{(i)}(\bx',\bth')}^*\right] &\approx
e^{i k (\bx \cdot \bth - \bx' \cdot \bth')}\notag\\
&\hspace{-0.2in} \times\exp{\Big( -\frac{3|\bm{P}_{{\bvth}} \tbx|^2 - 3
      |\bx-\bx^{(i)}(\bth)| \tbx \cdot \bm{P}_{\bvth}\tth +
      |\bx-\bx^{(i)}(\bth)|^{2} |\bm{P}_{{\bvth}}\tth|^2}{2
      X_{d,1}^2}\Big)},
\label{eq:MOMDi}
\end{align}
with $X_{d,1}$ defined as in \eqref{eq:defXd}, the notation \[\tbx =
\bx-\bx', \qquad \tth = \bth - \bth',
\] and the orthogonal projection 
\begin{equation}
\bm{P}_{\bvth} = I - \bvth \bvth^T.
\label{eq:orthbvth}
\end{equation} 

The second moments of the waves impinging on the scatterer at $\by$ are 
\begin{align}
\EE\left[u_1^{(i)}(\by,\bth) {u_1^{(i)}(\by,\bth')}^*\right] \approx
e^{i k \by \cdot \tth-\frac{|\bm{P}_{\bvth}\tth|^2}{2 \Theta_d^2}},
\label{eq:MOMDiy}
\end{align}
with dimensionless scale 
\begin{equation}
\Theta_d = \frac{X_{d,1}}{|\by-\by^{(i)}(\bvth)|} =
O\left(\frac{\la \sqrt{\ell}}{\sigma \sqrt{L^3}}\right) \ll \frac{\ell}{L} \ll 1
\label{eq:defThed}
\end{equation}
that defines the direction offset over
which the incoming  waves remain statistically correlated when
they reach the scatterer at $\by$.  By 
equations \eqref{eq:APA2} and  \eqref{eq:as1}, the scale $\Theta_d$ is much smaller
than the opening angle $\alpha$ of the cone $C$, 
\[
\Theta_d \ll \frac{\ell}{L} \ll \alpha
= O\left(\frac{a}{L} \right),
\]
which  is essential for obtaining statistically stable CINT
images, as we show later.

We state one more wave decorrelation result needed in the next sections. It says  
that the Green's function from the scatterer to the array and the  wave impinging on the scatterer are statistically decorrelated. This is because these waves traverse different parts of the random medium. More precisely,
let $\bx$ be a point in $A$ and $\bth$ a unit vector in the illumination cone $C$. We have 
\begin{align}
\EE\left[ G(\bx,\by;j \om) e^{i j k \gamma(\by,\bth)}
  \right] &\approx \EE\left[ G(\bx,\by;j \om)\right]\EE\left[e^{i j k \gamma(\by,\bth)}\right] \nonumber \\
  &\approx G_0(\bx,\by;j
\om)e^{-\frac{|\bx-\by|}{\ell_j^s} -
  \frac{|\by-\by^{(i)}(\bth)|}{\ell_j^s}} \approx 0. 
\label{eq:Mig2}
\end{align}

\subsection{Linearized data model in a random medium}
\label{sect:LINmod}
Using the weak nonlinearity assumption, we can write the 
solutions of equations \eqref{eq:waveu}--\eqref{eq:wavev} approximately as 
\begin{align}
u_1(\bx;\bth) &\approx u_1^{(i)}(\bx;\bth) + k^2 \left< \eta_1 \right> G(\bx,\by;\om)
u_1^{(i)}(\by;\bth), \label{eq:Born1} \\
u_2(\bx;\bth) &\approx (2 k)^2 \left< \eta_2 \right> G(\bx,\by;2\om)
[u_1^{(i)}(\by;\bth)]^2, \label{eq:Born2}
\end{align}
where we have represented the small scatterer by the net
susceptibilities $\left< \eta_j\right>$, given by the integral of $\eta_j$ over its small support, contained inside a ball centered at $\by$, 
of radius much smaller than $\la$,
\begin{equation*}
\left<\eta_j\right> = \int_{V} d \by \, \eta_j(\by) \ , \quad j = 1, 2.
\end{equation*}
The Green's function in equations \eqref{eq:Born1}--\eqref{eq:Born2}
propagates the waves in the medium, from the scatterer
to the array, and it is given by \eqref{eq:A5}. The direct wave
$u_1^{(i)}$ is the plane wave distorted by the random
medium, as given in equation \eqref{eq:A7}.

The random model of the data \eqref{eq:data1rand}-\eqref{eq:data2rand}
at the array is
\begin{align}
d_1(\bx,\bth) &= e^{i k \bx \cdot \bth} \left[ e^{i k
    \gamma(\bx,\bth)}-1\right] + k^2 \left< \eta_1 \right> G(\bx,\by;\om)
e^{i k \by \cdot \bth + ik
  \gamma(\by,\bth)}, \label{eq:dat1}\\ d_2(\bx,\bth) &= 4k^2
\left< \eta_2 \right> G(\bx,\by;2\om) e^{i 2k \by \cdot \bth + i 2 k
  \gamma(\by,\bth)}, \label{eq:dat2}
\end{align}
for $\bx \in \mathbb{A}$ and $\bth$ in the cone $C$ with axis along $\bvth$ and
opening angle $\alpha$. We neglect for simplicity the additive,
uncorrelated noise, which is much easier to handle than the random
medium distortions. 

\subsection{Analysis of migration imaging}
\label{sect:anal.2}
We assume in this and the following section that the number $N_{\bx}$ of sensors in the
array aperture $A$ is large, so that we can replace the sums over the
detectors by integrals over the aperture
\[
\sum_{s=1}^{N_{\bx}} \sim \int_{A}d \bx_\perp,
\]
where the symbol ``$\sim$'' denotes approximate, up to multiplication
by a constant. Recall that $\bx_\perp$ is the two dimensional vector 
in the square aperture $A$ of side $a$. 

We also approximate the sums over the incident
directions $\bth_q$ by integrals over the unit vectors $\bth$ in the
cone $C$, parametrized by the polar angle $\varphi \in (0, \alpha)$ between $\bth$ and $\bvth$ and the 
azimuthal angle $\beta \in [0,2 \pi]$ 
\[
\sum_{q=1}^{N_{\theta}} \sim \int_{C} d \bth = \int_0^\alpha d \varphi \, \sin \varphi \int_0^{2 \pi} d \beta.
\]
The migration imaging function  \eqref{eq:LS4}  at search points $\bys \in R$ is modeled (up to
multiplicative constants) by
\begin{align}
\IMj(\bys) = \int_{A} d \bx_\perp \int_{C} d \bth \, 
d_j(\bx,\bth) {G_0^\star(\bys,\bx;j\om)} e^{-i j k \by \cdot \bth},
\label{eq:Mig1}
\end{align}
with $d_j$ given by \eqref{eq:dat1}-\eqref{eq:dat2}.  We describe first its focusing in homogeneous media, 
and then consider  random media.

\subsubsection{Homogeneous media}
When the wave speed equals the constant $c_0$,  the Green's function equals $G_0$ and there is no distortion of the direct wave, so the 
first term in \eqref{eq:dat1} cancels out. Substituting the resulting data model 
in \eqref{eq:Mig1}, we obtain the migration imaging function  
\begin{align}
\HMj(\bys) = (jk)^2 \left< \eta_j \right> \int_{A} d \bx_\perp \, 
G_0(\bx,\by;j \om) {G_0^\star(\bys,\bx;j\om)} \int_{C} d \bth \, e^{i j k \bth \cdot (\by-\bys)}.
\label{eq:RH0}
\end{align}
It has a separable form, given by the product of two integrals over the array aperture $A$ and the cone 
$C$ of illuminations. 

The integral over the aperture is 
\begin{align}
\int_{A} d \bx_\perp \, 
G_0(\bx,\by;j \om) {G_0^\star(\bys,\bx;j\om)} &\approx \frac{a^2}{L^2} e^{ i j k \Big( y_\parallel - y^R_\parallel + 
\frac{|\by_\perp|^2 - |\bys|^2}{2L} \Big) } \notag\\
&\times \mbox{sinc}\Big(\frac{j k a(y-y^R)_1}{2L}\Big)\mbox{sinc}\Big(\frac{j k a(y-y^R)_2}{2L}\Big),
\label{eq:RH1}
\end{align}
where we used the paraxial approximation \eqref{eq:APAR3}, and indexed by $1$ and $2$ 
the components of $\by_\perp$ and $\bys_\perp$ in the plane orthogonal to $\bn$,
for coordinate axes parallel to the sides of the square aperture. 
This is the classic calculation of the point spread function of time reversal in homogeneous media.
It localizes the scatterer in the plane orthogonal to $\bn$ with resolution of the order $\la L/(ja)$.

The integral over the cone $C$ of illuminations is 
\begin{align}
\int_{C} d \bth \, e^{i j k \bth \cdot (\by-\bys)} &= 
\int_0^\alpha d \varphi \, \sin \varphi \, e^{i j k \cos \varphi \, \bvth \cdot (\by-\bys)}\int_0^{2 \pi} d \beta \, e^{i j k  
\sin \varphi \cos \beta |\bm{P}_{\bvth}(\by-\bys)|} ,
\label{eq:RH2}
\end{align}
and we can simplify it using the assumptions \eqref{eq:APA2} and  \eqref{eq:APAR} on the small opening angle $\alpha$ and the linear size $r$ of the search domain. We approximate 
\[
k \cos \varphi \, \bvth \cdot (\by-\bys) = k \bvth \cdot (\by-\bys) + O\left(\frac{r a^2}{\la L^2}\right) \approx k \bvth \cdot (\by-\bys) ,
\]
and 
\[
k   \sin \varphi \, \cos \beta |\bm{P}_{\bvth}(\by-\bys)| =k \varphi \cos \beta |\bm{P}_{\bvth}(\by-\bys)| + O \left(\frac{a^3 r}{\la L^3} \right) \approx 
k \varphi \cos \beta |\bm{P}_{\bvth}(\by-\bys)| ,
\]
and  obtain that 
\begin{align}
\int_{C} d \bth \, e^{i j k \bth \cdot (\by-\bys)} & \approx e^{i j k \bvth \cdot (\by-\bys)}
\int_0^\alpha d \varphi \, \varphi \, \int_0^{2 \pi} d \beta \, e^{i j k  
\varphi \cos \beta |\bm{P}_{\bvth}(\by-\bys)|} \nonumber \\
 &= 2 \pi \,{e^{i j k \bvth \cdot (\by-\bys)}} \int_0^\alpha d \varphi \, \varphi J_0\left(j k \varphi |\bm{P}_{\bvth}(\by-\bys)|\right) \nonumber \\
 & = 2 \pi \,{e^{i j k \bvth \cdot (\by-\bys)}} \alpha \frac{J_1 \left(j k \alpha |\bm{P}_{\bvth}(\by-\bys)|\right)}{j k |\bm{P}_{\bvth}(\by-\bys)|},
\label{eq:RH3}
\end{align}
where $J_q$ are the Bessel functions of the first kind for $q = 0, 1$.  This expression is large when
\[
|\bm{P}_{\bvth}(\by-\bys)| = O \left(\frac{\la}{j \alpha}\right) = O\left(\frac{\la L}{j a}\right),
\]
and gives  the focusing in the plane orthogonal to $\bvth$, and therefore along $\bn$.

The imaging function is the product of \eqref{eq:RH1} and \eqref{eq:RH2} and focuses at $\by$ with resolution $\la L/(ja)$.

\subsubsection{Migration image in random media}
To analyze the behavior of the migration imaging function in random media, we calculate in appendix \ref{ap:mig} its expectation and standard deviation. The result is summarized as follows.
The expectation $\EE\big[\IMj(\by)\big]$ of the migration imaging function \eqref{eq:Mig1} evaluated at the 
scatterer location $\by$ is much smaller than its standard deviation. {The signal to noise ratio (SNR), which is the ratio of the expectation to the standard deviation, satisfies 
\begin{equation}
{\rm SNR}\big[\IMj(\by)\big] \ll 1.
\label{eq:expsmall}
\end{equation}
This result means that we cannot draw any conclusion about the focusing of the migration image 
by studying its statistical expectation. Because the waves that reach the array are
randomized in our scaling, as stated in section
\ref{sect:anal.1.1}, the migration image is also randomized, and has very large fluctuations with respect to its mean. 
This is manifested in practice by the fact that the image may not be focused and reproducible, since it changes unpredictably 
 with the realizations of the medium. 
There is no mechanism for mitigating the wave randomization 
in migration. The integration over the 
array aperture and over the illumination directions only takes care 
of additive and  uncorrelated noise, but it cannot deal with the large random wave distortions due to scattering in the medium.
This is why migration is not statistically stable. 

\subsection{Analysis of CINT imaging}
\label{sect:anal.3}
The CINT imaging function is defined in \eqref{eq:CINT6}, with the sums replaced by integrals over the  aperture $A$
and the cone $C$, except for one modification, namely that the waves decorrelate over offsets $\tth = \bth-\bth'$ in the plane orthogonal  to $\bvth$, so we replace  in \eqref{eq:CINT6} 
\[
\Phi \Big(\frac{\tth}{\Theta}\Big) \to \Phi \Big(\frac{\bm{P}_{\bvth} \tth}{\Theta}\Big),
\]
where $\bm{P}_{\bvth}$ is as defined in \eqref{eq:orthbvth}.
We also take a Gaussian window $\Phi(\bm{z}) = \exp(-{|\bm{z}|^2}/{2})$ to simplify the calculations. The thresholding parameters 
$X$ and $\Theta$ are of the same order as the decoherence scales 
$X_{d,1}$ and  $\Theta_{d}.$

The CINT image is formed with cross-correlations of the measurements at points 
$
\bx \pm {\tbx}/{2} \in \mathbb{A},$ and for incident  plane waves with unit wave vectors
$
\bth \pm \tth/{2} \in C.
$
In our system of coordinates we have $\bx = (\bx_\perp, 0)$ and 
$\tbx = (\tbx_\perp,0)$, with 
\begin{equation}
(\bx_\perp,\tbx_\perp) \in  \left\{ (\bm{z}, \tilde{\bm{z}}) \in \mathbb{R}^4 : \ |z_j| \le \frac{a}{2}, 
~ \mbox{and} ~ ~|\tilde z_j| \le \min\{a - 2|z_j|, 3 X \}, ~ j = 1, 2 \right\}.
\label{eq:defSetA}
\end{equation}
Here we used the fact that the offset is limited by the essential support of the Gaussian window 
\[
\Phi(\tbx) = \exp\Big(-\frac{|\tbx|^2}{2 X^2} \Big) =  \exp \Big(-\frac{|\tbx_\perp|^2}{2 X^2} \Big),
\] 
which is three times its standard deviation. Note that since $X \sim X_{d,1} \ll a$,  the offsets  $\tbx_\perp$ are limited by $3 X$ 
for most center points $\bx_\perp$, so we can obtain a good approximation of the imaging function by using the simpler set 
\begin{equation}
\cA = \left\{ (\bm{z}, \tilde{\bm{z}}) _\xi\in \mathbb{R}^4 : \ |z_j| \le \frac{a}{2}, 
~ \mbox{and} ~ ~|\tilde z_j| \le 3 X, ~ j = 1, 2 \right\}.
\label{eq:defSetAS}
\end{equation}
We denote by  $\iint_{\cA} d \bx_\perp d \tbx_\perp$ the integral over $\cA$.

To define the set that supports $\bth$ and $\tth$, we use the orthonormal basis $\{\bvth, \bet, \bzet\}$, with 
vector $\bet$ aligned with the projection $\bm{P}_{\bvth} \bth$, so that 
\begin{equation}
\bth = |\bth| (\bvth \cos \varphi  + \bet \, \sin \varphi ), \quad \varphi \in (0,\alpha).
\label{eq:Th1}
\end{equation}
Since $|\bth \pm \tth/2| = 1$, we  have 
\begin{equation}
\bth \cdot \tth = 0, \quad |\bth| = \sqrt{1 - |\tth|^2/4},
\end{equation}
and using the  decomposition
\begin{equation}
\tth = \tilde{\theta}_{\vartheta} \bvth + \tilde{\theta}_\xi \bet + \tilde{\theta}_\zeta \bzet,
\label{eq:Th2}
\end{equation}
we can solve for the component of $\tth$ along the axis of the cone $C$,
\begin{equation}
\tilde{\theta}_{\vartheta} = - \tan \varphi \, \tilde \theta_\xi.
\label{eq:Th3}
\end{equation} 
This yields 
\begin{equation}
|\tth|^2 = \tilde \theta_{\zeta}^2 + \frac{\tilde \theta_{\xi}^2}{\cos^2 \varphi} = \tilde \theta_{\zeta}^2  + \tilde \theta_{\xi}^2  + O(\alpha^2 \Theta_d^2), 
\quad |\bth| = 1 + O(\Theta_d^2), 
\label{eq:Th4}
\end{equation} 
where we used
\[
|\bm{P}_{\bvth} \tilde \theta|^2 =  \tilde \theta_{\zeta}^2  + \tilde \theta_{\xi}^2 = O( \Theta_d^2),
\]
due to the Gaussian thresholding window with $\Theta = O(\Theta_d)$. We write then that $(\bth, \tth) \in \cC$, 
the set defined by vectors $\bth$ of the form \eqref{eq:Th1}, with norm as in \eqref{eq:Th4}, and 
\begin{equation}
\tth = \tilde \theta_\xi \left( \bet - \bvth \tan \varphi \right) + \tilde \theta_\zeta \bzet.
 \label{eq:Th5}
 \end{equation} 
We parametrize $\cC$ by  the polar angle $\varphi \in (0, \alpha)$, the azimuthal angle $\beta \in (0, 2 \pi)$, and  the components $\tilde \theta_\xi$ 
and $\tilde \theta_\zeta$ of $\tth$. The angle  $\beta$ determines the unit vectors $\bet$ and $\bzet$ in polar coordinates, in the 
plane orthogonal to $\bvth$. 
We also denote by  $\iint_{\cC} d \bth d \tth$ the integral over $\cC$.

The analysis of CINT in homogeneous media is not interesting. This is because  the windowing of the detector and direction offsets 
is not necessary, and once we remove it, the CINT imaging function becomes the square of the migration function.  
We analyze separately the imaging of the linear and quadratic susceptibilities in random media.  The calculations are similar, except that 
in the linear case  data \eqref{eq:dat1} have an extra term due to the randomly distorted direct wave. We begin   with the imaging of the 
quadratic susceptibility, which uses the simpler data  model \eqref{eq:dat2}.
\subsubsection{Imaging of the quadratic susceptibility}
\label{sect:QCint}
The model of \eqref{eq:CINT1} at imaging point $\bys$ and the second harmonic frequency $2 \om$ is 
\begin{align}
b_2\Big(\bx\pm \frac{\tbx}{2},\bth\pm \frac{\tth}{2},\bys\Big) = 4 k^2 \left< \eta_2 \right> G\Big(\bx\pm \frac{\tbx}{2},\by;2 \om\Big) 
G_0^\star\Big(\bx\pm \frac{\tbx}{2},\bys;2 \om \Big) \nonumber \\
\times \exp \Big[i 2 k \Big(\bth + \frac{\tth}{2}\Big) \cdot (\by-\bys) + i 2 k \gamma\Big(\by,\bth\pm\frac{\tth}{2}\Big)\Big], 
\label{eq:CI1}
\end{align}
for $\bx = (\bx_\perp,0)$ and  $\tbx = (\tbx_\perp, 0)$ with  $(\bx_\perp, \tbx_\perp) \in \cA$, and  $(\bth, \tth) \in \cC.$
The image is formed with the cross-correlations of \eqref{eq:CI1}
\begin{align}
\QC(\bys) = \iint_{\cA} d \bx_\perp d\tbx_\perp \iint_{\cC} d \bth d \tth \, \exp \Big(-\frac{|\tbx_\perp|^2}{2 X^2}-\frac{|\bm{P}_{\bvth} \tth|^2}{2 \Theta^2} \Big) \nonumber \\
\times 
b_2 \Big(\bx+ \frac{\tbx}{2},\bth+ \frac{\tth}{2}, \bys \Big)
b_2^\star  \Big(\bx -\frac{\tbx}{2}, \bth - \frac{\tth}{2}, \bys\Big).
\label{eq:CI3}
\end{align}
Its focusing and statistical stability are described by the following results derived in appendix \ref{ap:CINT}:

The expectation of \eqref{eq:CI3} is given by 
\begin{align}
\EE\left[\QC(\bys)\right] &\approx \frac{(2 \pi)^3}{2} \left(4 k^2 \left< \eta_2 \right>^2  \alpha \Theta_e \frac{a X_e}{L^2} \right)^2\notag\\
&\times\exp\left[-\frac{1}{2} \left(\frac{2 k X_e |\by_\perp-\bys_\perp|}{L}\right)^2 - \frac{1}{2} 
\left(2 k \Theta_e |\bm{P}_{\bvth}(\by-\bys)|\right)^2\right],
\label{eq:CI4}
\end{align}
where $X_e$ and $\Theta_e$ are defined by 
\[
\frac{1}{X_e^2} = \frac{1}{X^2} + \frac{1}{X_{d,2}^2}, \quad \frac{1}{\Theta_e^2} = \frac{1}{\Theta^2} + \frac{4}{\Theta_d^2}.
\]
The SNR of  \eqref{eq:CI3} evaluated at the scatterer location is 
\begin{equation}
SNR\left[\QC(\by)\right] \sim (a/\ell)^2.
\label{eq:SNRCint2}
\end{equation}
Since in our scaling $\ell \ll a$, the SNR   is high, meaning that 
$
\QC(\bys) \approx \EE\left[\QC(\bys)\right],
$
for $\bys$ near $\by$. This is the statement of statistical stability. The focusing of the image is 
determined by the exponential in \eqref{eq:CI4}. The first term gives the focusing in the plane of the array,
and the second in the plane orthogonal to $\bvth$, where we recall that $\bvth \perp \bn$, and 
$\bn$ is normal to the array. By assumption \eqref{eq:decorestim} we have  $X_e \sim X_{d,2}$ and $\Theta \sim \Theta_d$, 
and using definition \eqref{eq:defThed} of $\Theta_d$ we conclude that the CINT resolution is
\begin{equation}
|\by-\bys| \le O\left(\frac{\la L}{X_e} \right) \gg \frac{\la L}{a}.
\label{eq:CI5}
\end{equation}
Similar to the results in \cite{borcea2011enhanced,borcea2006adaptive}, we conclude that  the cost of statistical stability comes at the expense of 
resolution\footnote{Note that the peak of the image can be observed in the search region $R$ with linear size $r$ satisfying \eqref{eq:APAR}, because 
\[
\frac{\la L/X_{e}}{\la L^2/a^2} = O\left(\frac{a^2}{L X_{d,2}} \right) = O\left(\frac{a^2 \sigma}{\la \sqrt{\ell L}}\right) \ll \frac{a^2}{\sqrt{\la L^3}} \ll 1.
\]
}
which is lower than in  homogeneous media. 

\subsubsection{Imaging of the linear susceptibility}
\label{sect:LCint}
The model of \eqref{eq:CINT1} at imaging point $\bys$ and at frequency $\om$ is given by 
\begin{align}
b_1\Big(\bx\pm \frac{\tbx}{2},\bth\pm \frac{\tth}{2},\bys\Big) = \Big[e^{i k \gamma\big(\bx\pm \frac{\tbx}{2},\bth\pm \frac{\tth}{2}\big)}-1\Big] G_0^\star\Big(\bx\pm \frac{\tbx}{2},\bys;\om\Big)e^{i k \big(\bth \pm \frac{\tth}{2}\big)\cdot \big(\bx \pm \frac{\tbx}{2} - \bys\big) }\nonumber \\
+ k^2 \left< \eta_1 \right> G\Big( \bx \pm \frac{\tbx}{2},\by;\om \Big) G_0^\star \Big( \bx \pm \frac{\tbx}{2},\bys;\om\Big) e^{i k \big(\bth \pm \frac{\tth}{2}\big) \cdot(\by-\bys) + i k \gamma\big(\by,\bth+\frac{\tth}{2}\big)}.
\label{eq:CI6}
\end{align}
The first term is due to the uncompensated direct wave which has not interacted with the scatterer at $\by$. The second term is the 
useful one in inversion. 
The imaging function is given by the superposition of the cross-correlations of \eqref{eq:CI6},
\begin{align}
\LC(\bys) &= \iint_{\cA} d \tbx d \tbx_\perp \iint_{\cC} d \bth d \tth \, e^{-\frac{|\tbx_\perp|^2}{2 X^2}-\frac{|\bm{P}_{\bvth}\tth|^2}{2 \Theta^2}}\notag\\
&\times
b_1 \Big(\bx+ \frac{\tbx}{2},\bth+ \frac{\tth}{2}, \bys \Big)
b_1^\star  \Big(\bx -\frac{\tbx}{2}, \bth - \frac{\tth}{2}, \bys\Big),
\label{eq:CI8}
\end{align}
with the same notation as in the previous section. 

We derive in appendix \ref{ap:CINT} the expression of the 
expectation of \eqref{eq:CI8},
\begin{align}
\EE\Big[\LC(\bys)\Big] &\approx \frac{(2 \pi)^3}{2} \left(k^2 \left< \eta_1 \right>^2 \alpha \Theta_e \frac{a X_e}{L^2} \right)^2\notag\\
&\times
\exp \left[-\frac{1}{2} \left(\frac{k X_e |\by_\perp-\bys_\perp|}{L}\right)^2 - \frac{1}{2} \left(k \Theta_e
|\bm{P}_{\bvth}(\by-\bys)|\right)^2\right],
\label{eq:CI9}
\end{align}
with $X_e$ and $\Theta_e$ defined by 
\[
\frac{1}{X_e^2} = \frac{1}{X^2} + \frac{1}{X_{d,1}^2}, \quad \frac{1}{\Theta_e^2} = \frac{1}{\Theta^2} + \frac{1}{\Theta_d^2}.
\]
This expression is like \eqref{eq:CI4}, except that in that equation we had the double frequency.
We also show in appendix \ref{ap:CINT} that  the SNR of \eqref{eq:CI8} evaluated at the scatterer location is  large, 
\begin{equation}
SNR\Big[\LC(\by)\Big] \sim \left(\frac{a}{\ell}\right)^2,
\end{equation}
meaning that the imaging function is statistically stable in the vicinity of the scatterer at $\by$.

\begin{remark}
{\rm The first term in \eqref{eq:CI6}, due to the uncompensated direct wave, does not play a role in these results because we limit the search points $\bys$ to the small region $R$ centered at $\by$. 
If we searched  in the whole domain,  we would see that due to  the direct wave, the expectation of the imaging function \eqref{eq:CI8} is large at points $\bys$ near the array. Moreover, the 
set of such points grows as we increase the aperture size and the opening angle of the cone of illuminations, 
in the sense that the larger these are, the further the points from the array that contribute to the image. 
We refer to appendix \ref{ap:CINTLin} for the justification of this statement. The result is of interest because it says that while in general it is advantageous to have a diverse set of incident directions and a larger aperture, this is not so when we cannot eliminate from the measurements the direct waves 
that have not interacted with the scatterers that we wish to image. These waves lead to spurious image peaks that cover a larger and larger
neighborhood of the array as we increase the opening angle of the cone of illuminations and the aperture, and make it difficult to 
locate the scatterers unless we know approximately where to search in a favorable position.}
\end{remark}

\section{Numerical results}
\label{sect:numerics}
\setcounter{equation}{0}
In this section we present numerical results of migration and CINT imaging using data calculated by solving the nonlinear equations \eqref{eq:waveu}--\eqref{eq:wavev},
as explained in appendix \ref{app:numerics}. To decrease the computational cost, the numerical study is performed in two dimensions. 
The setup is as shown in Figure \ref{imsetup}, but the scaling regime is different than that used in the analysis presented in section \ref{sect:anal}.
There are two reasons for this choice. The first is that the scaling necessary for analysis requires very long distances of propagation of the waves, over many wavelengths, which makes  the forward solver described in appendix \ref{app:numerics} prohibitively expensive computationally. The second reason is that we wish to explore  a different scattering regime, that is difficult to analyze theoretically, and yet gives qualitatively similar results to those predicted by the 
theory in section \ref{sect:anal}.

We consider a square domain $V$ of side length $20\lambda$, where $\lambda$ is the wavelength of the incident field. The array covers the entire bottom side of $V$, and the system 
of coordinates has its origin at the center of the array, with the $x_1$-axis pointing horizontally to the right. The domain $V$ contains two
small scatterers treated as disks of radius $0.1\lambda$. We indicate their true location in the images in Figures \ref{HK}-\ref{entire1} with black circles. The linear susceptibility of the scatterers 
is $1$, and the quadratic susceptibility is $0.01$. We take a wide cone  $C$ of incident directions, parametrized by the angle $\varphi \in [-\frac{\pi}{4},\frac{\pi}{4}]$, with center direction $\bvth$ pointing horizontally, to the right. Unless indicated otherwise in the caption of the figures, we use $20$ incident angles and $81$ sensors in the aperture $a = 20 \la$.

\begin{figure}[t]
\centering
\includegraphics[width=0.42\linewidth]{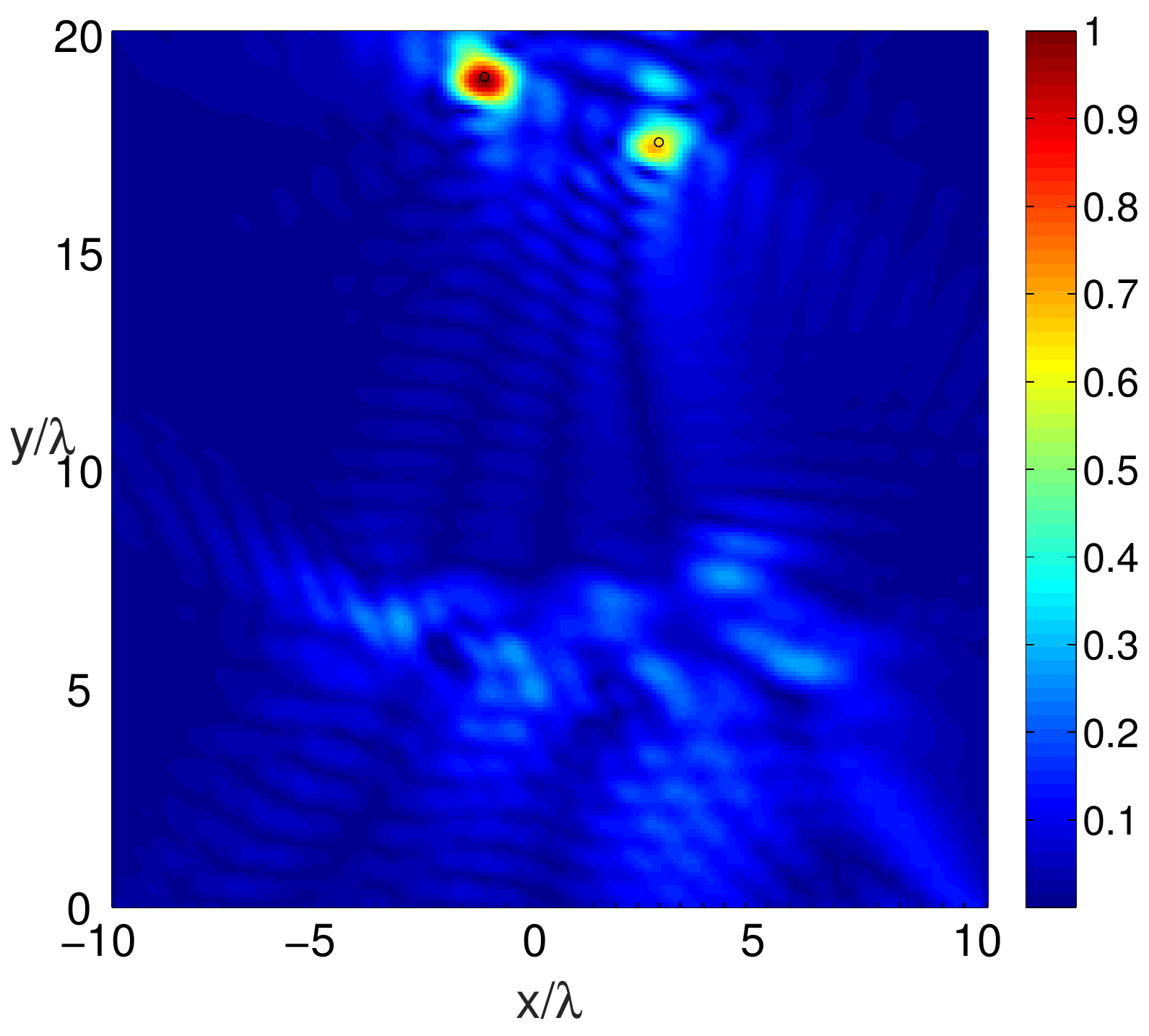}
 \includegraphics[width=0.42\linewidth]{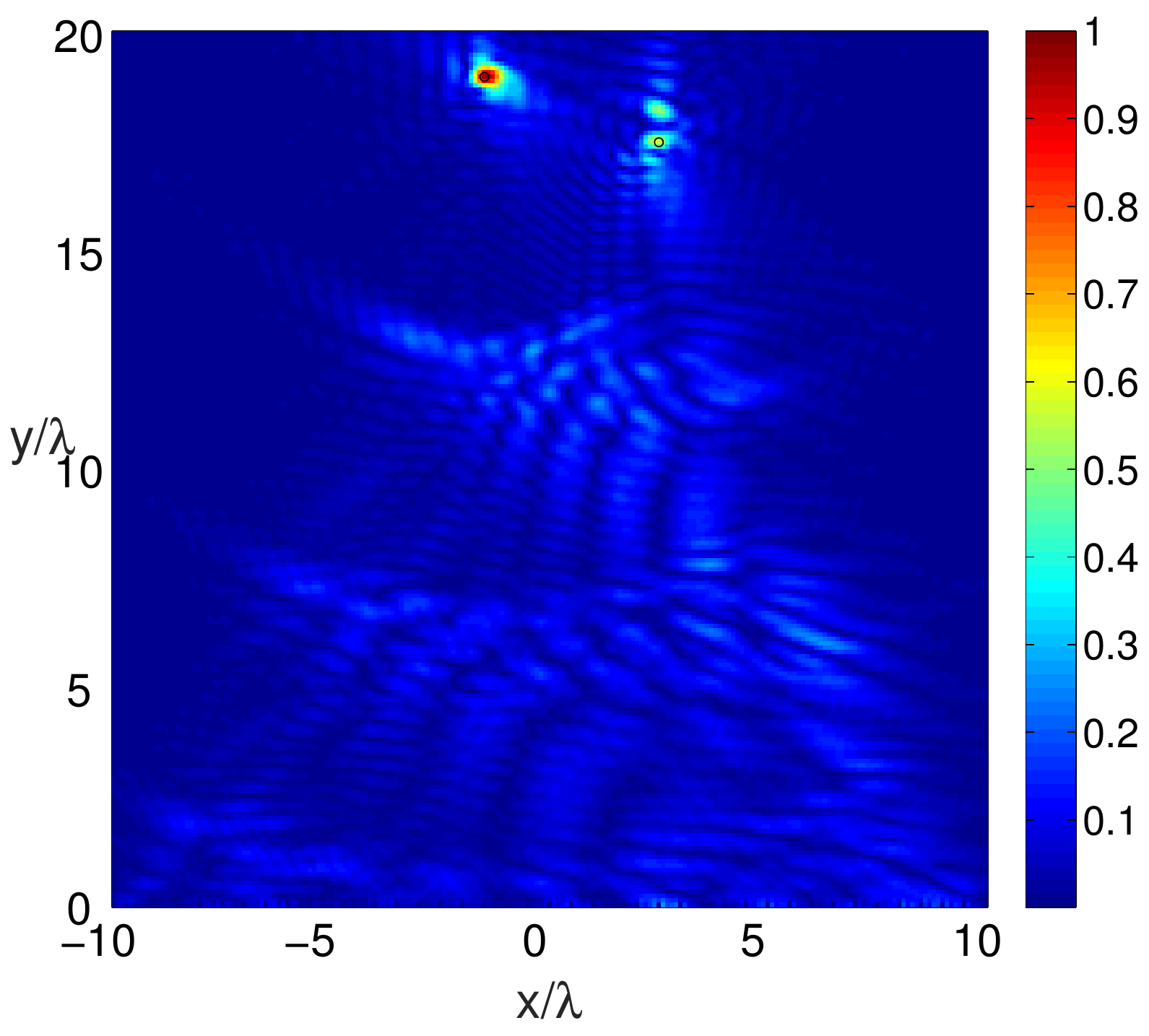}
\vspace{-0.1in}\caption{Migration images of $\eta^{(1)}$ (left) and $\eta^{(2)}$ (right) in a homogeneous medium.}
\label{HK}
\end{figure}
Migration images of $\eta_1$ and $\eta_2$ in a homogeneous medium are shown in Figure~\ref{HK}. They peak at the scatterer locations. We observe that the resolution of the images is of the order of the wavelength. Since this is smaller for the second harmonic, the resolution of the image of $\eta_2$ is better. 

\begin{figure}[t]
\centering
\includegraphics[width=0.45\linewidth]{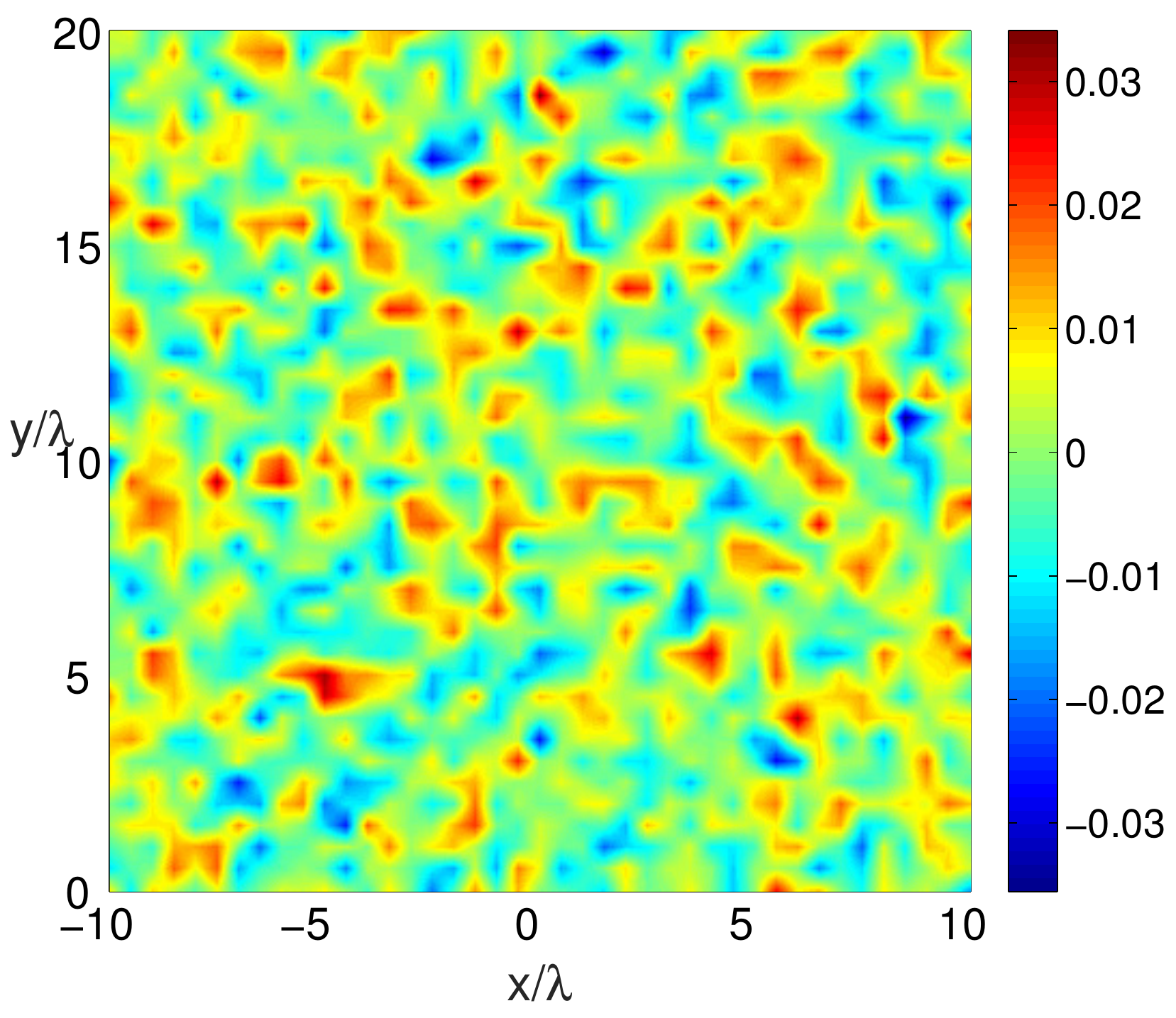}
\vspace{-0.1in}
\caption{One realization of the process $\mu$.
}
\label{back1}
\end{figure}

We study imaging in a random medium, generated numerically with random Fourier series \cite{devroye2006nonuniform} for the Gaussian autocorrelation function \eqref{eq:A2} with correlation length $\ell = 0.3 \la$, and standard deviation  $\sigma = 0.01\times(4 \pi)$.  We display in Figure 
\ref{back1} one realization of $\eta$. The migration and CINT images of $\eta_1$ and $\eta_2$ in a small search region near the scatterer are shown in Figures ~\ref{MC1} and \ref{MC2}, for two realizations of the random medium. The thresholding parameters in the CINT image formation are  $X_1=2X_2=7\lambda$ and $\Theta={\pi}/{5}$.
We observe that as predicted by the theory in section \ref{sect:anal}, the migration images change significantly from one realization to another, whereas  CINT images of $\eta_1$ and $\eta_2$  are almost the same and peak near the scatters.  The figures also show that  the CINT images are blurrier than those in homogeneous media displayed in Figure \ref{HK}.  
\begin{figure}[t]
\centering
\includegraphics[width=0.4\linewidth]{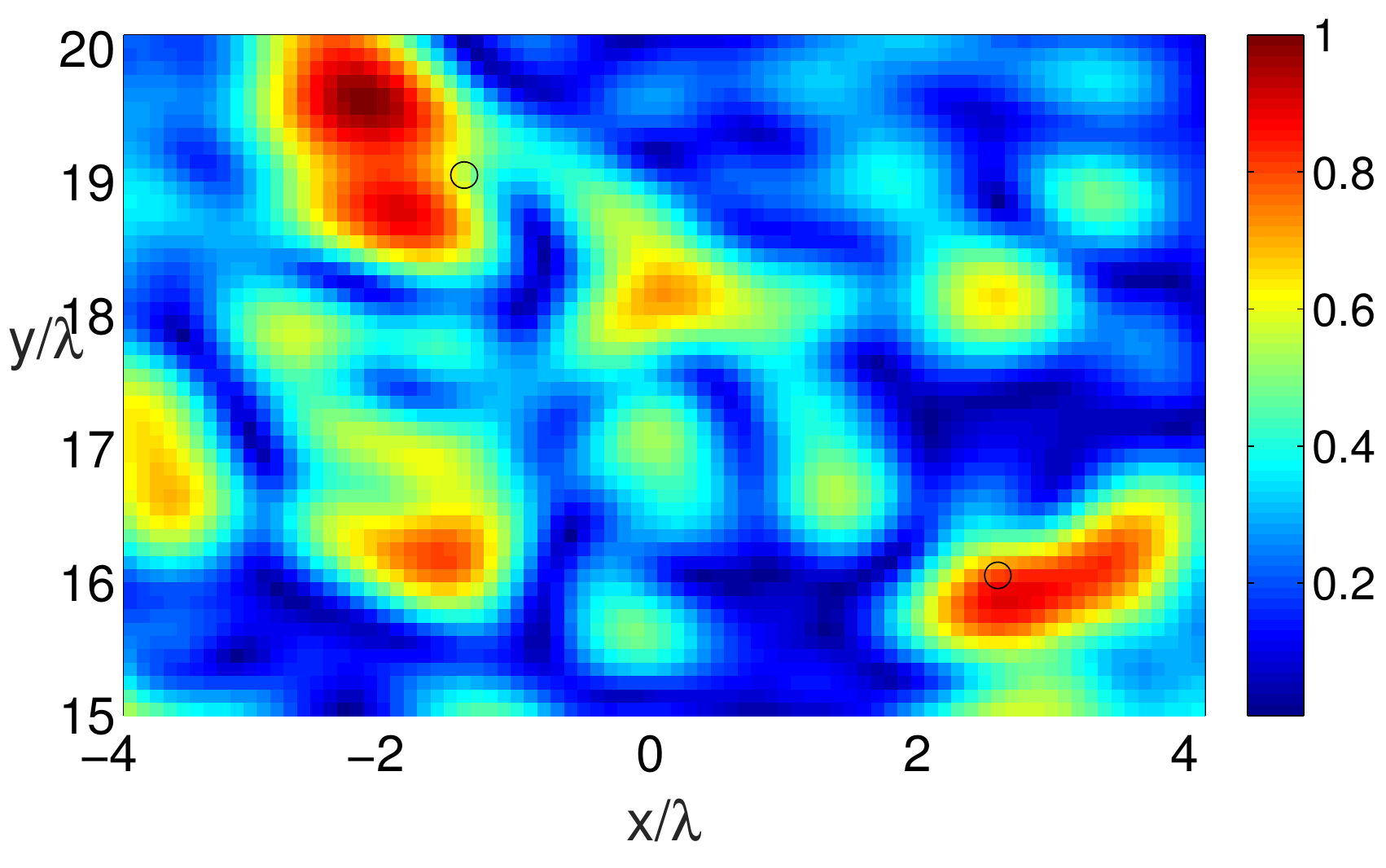} 
 \includegraphics[width=0.4\linewidth]{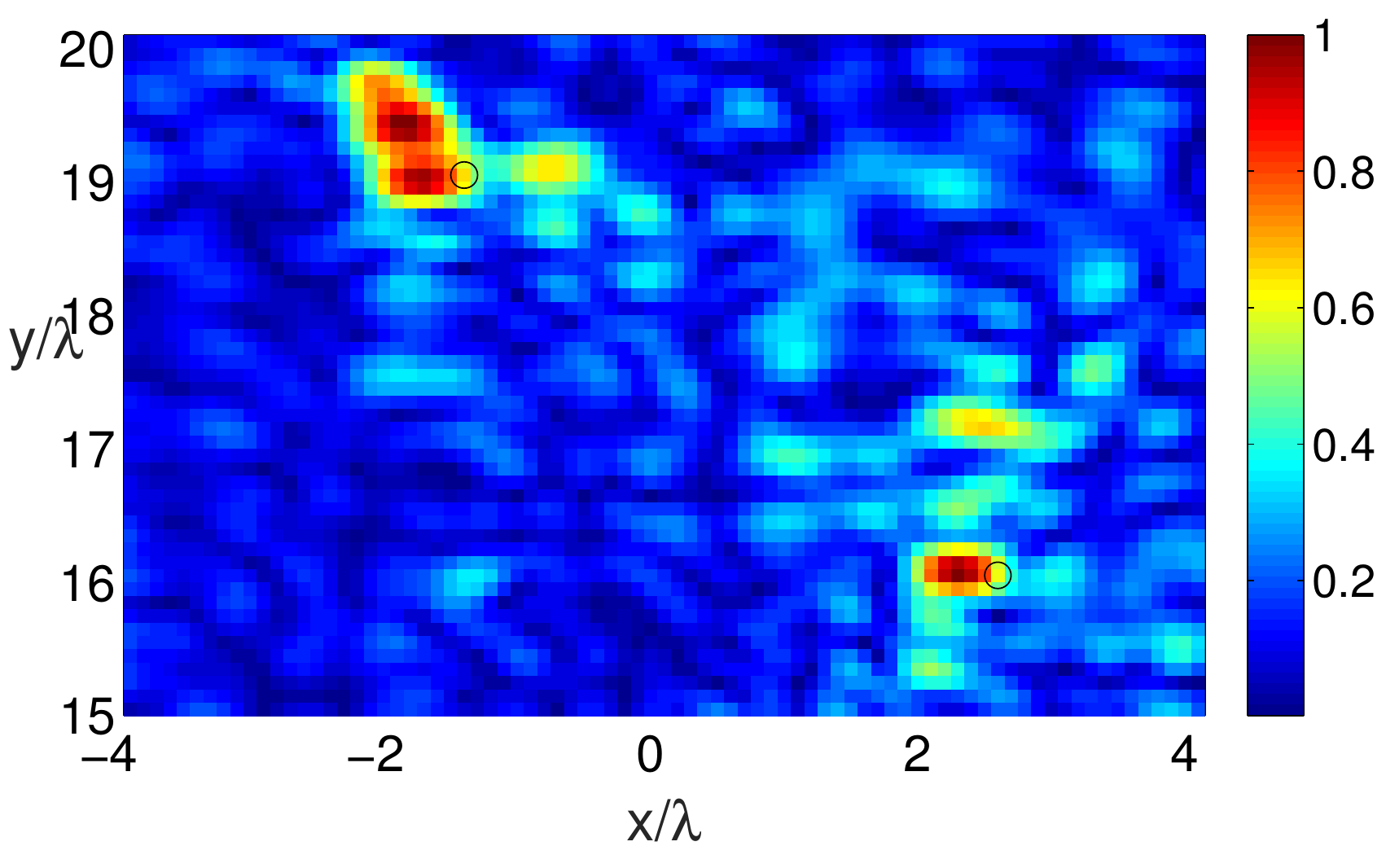}
 \\
 \vspace{-0.in}
\includegraphics[width=0.4\linewidth]{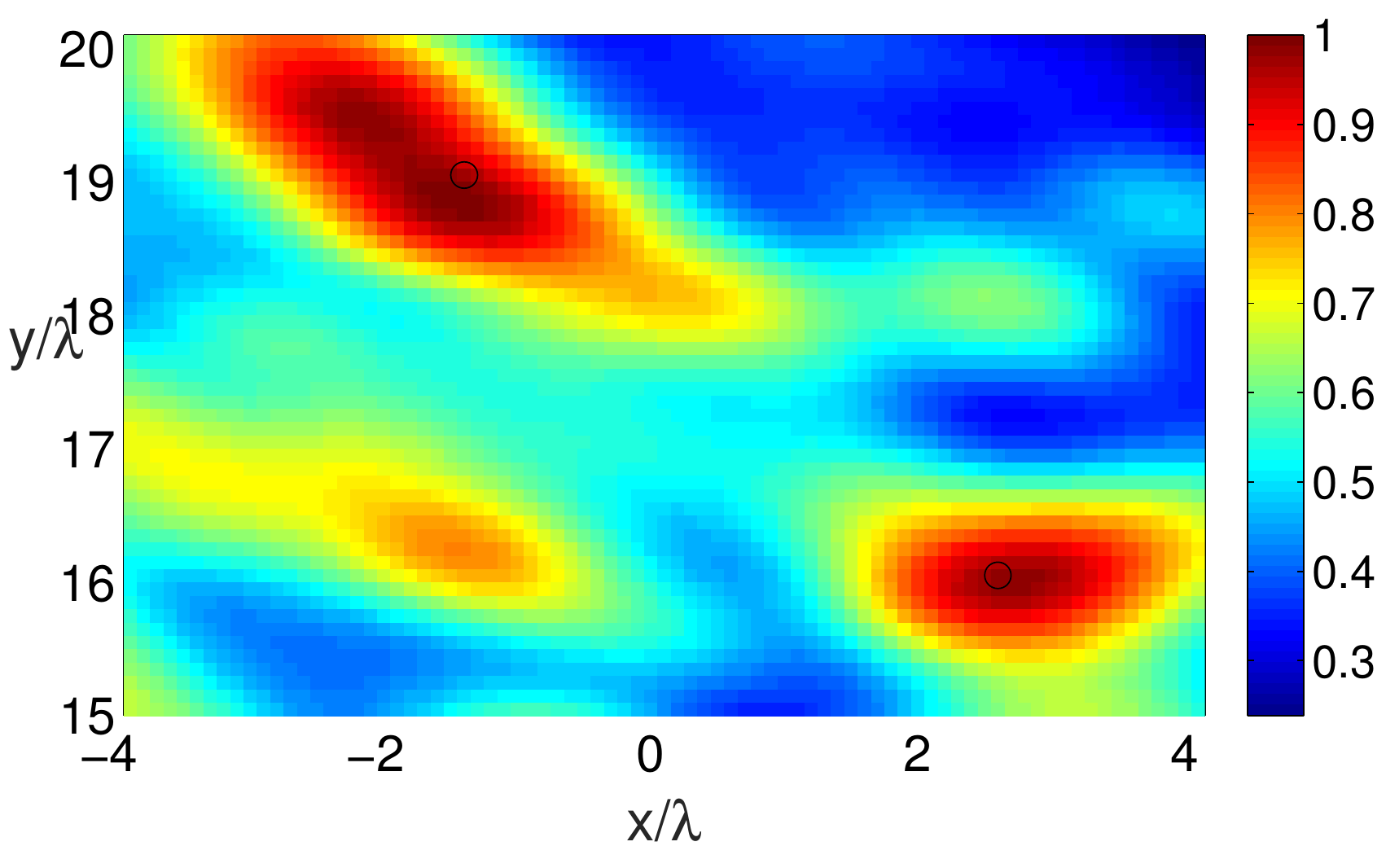}
\includegraphics[width=0.4\linewidth]{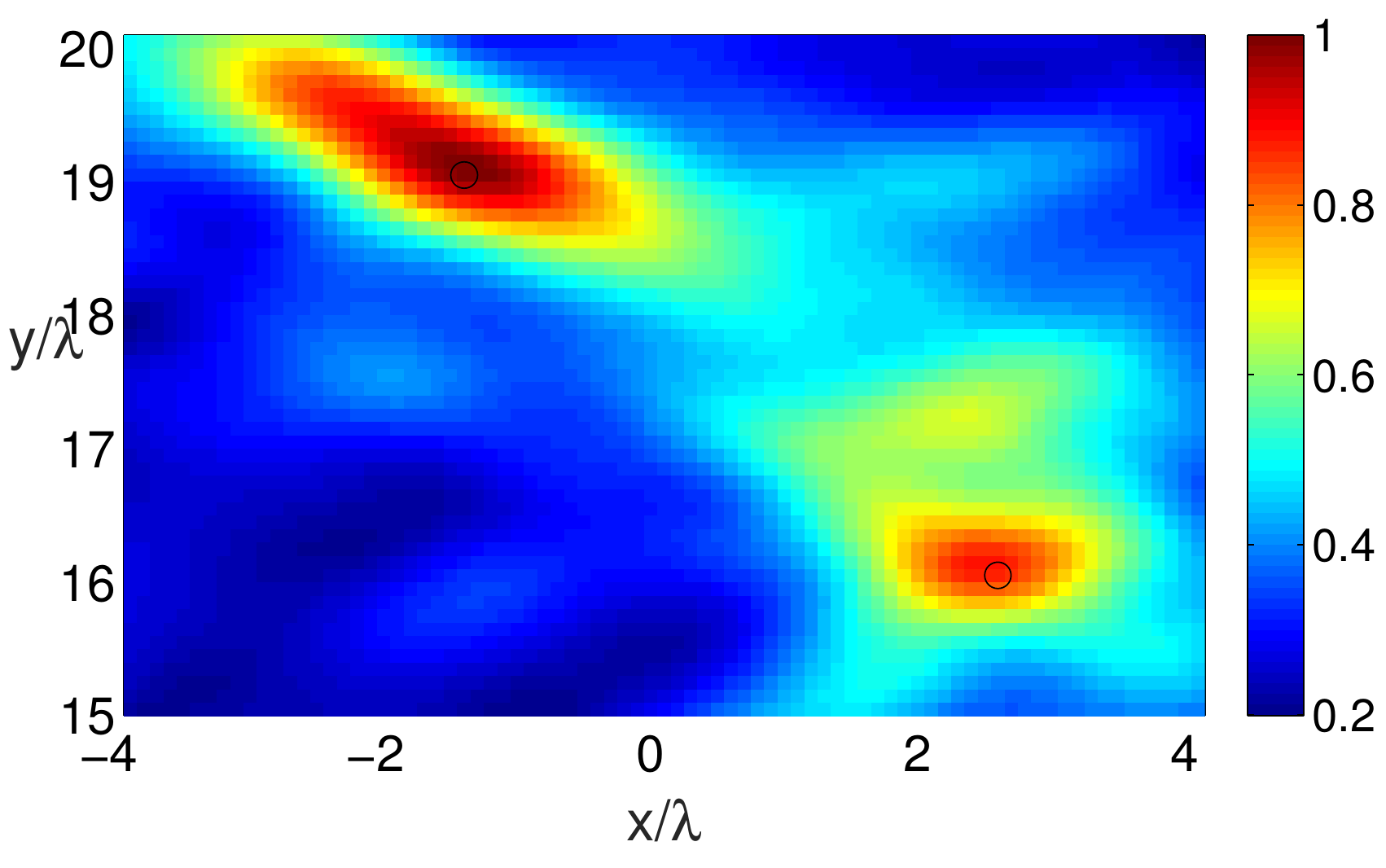}
\vspace{-0.1in}
\caption{The migration (top row) and CINT (bottom) row images of the linear susceptibility $\eta_1$ in two realizations of the random medium. The scatterer location is shown with a black circle.}
\label{MC1}
\end{figure}

\begin{figure}[t]
\centering
\centering
\includegraphics[width=0.4\linewidth]{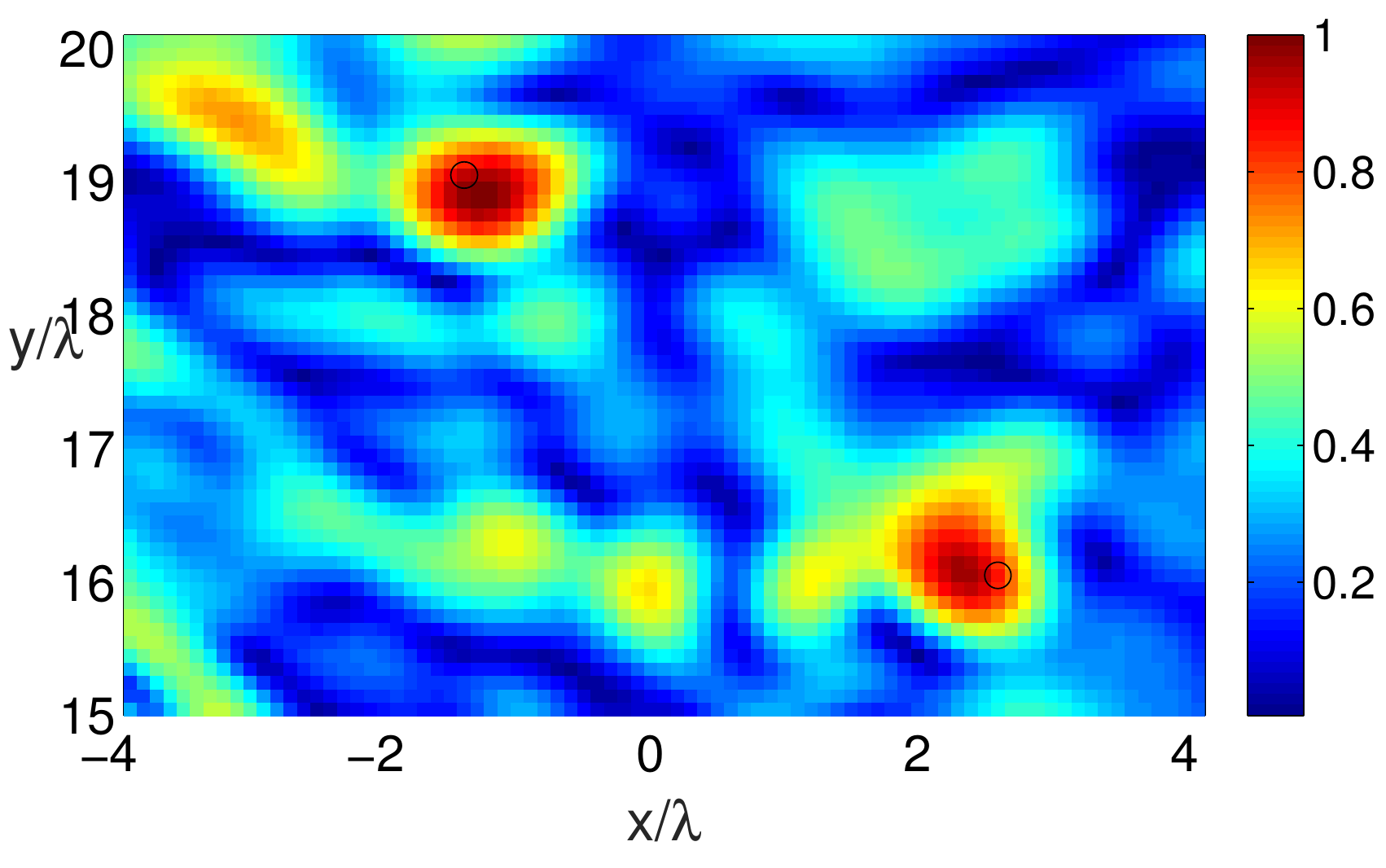} 
 \includegraphics[width=0.4\linewidth]{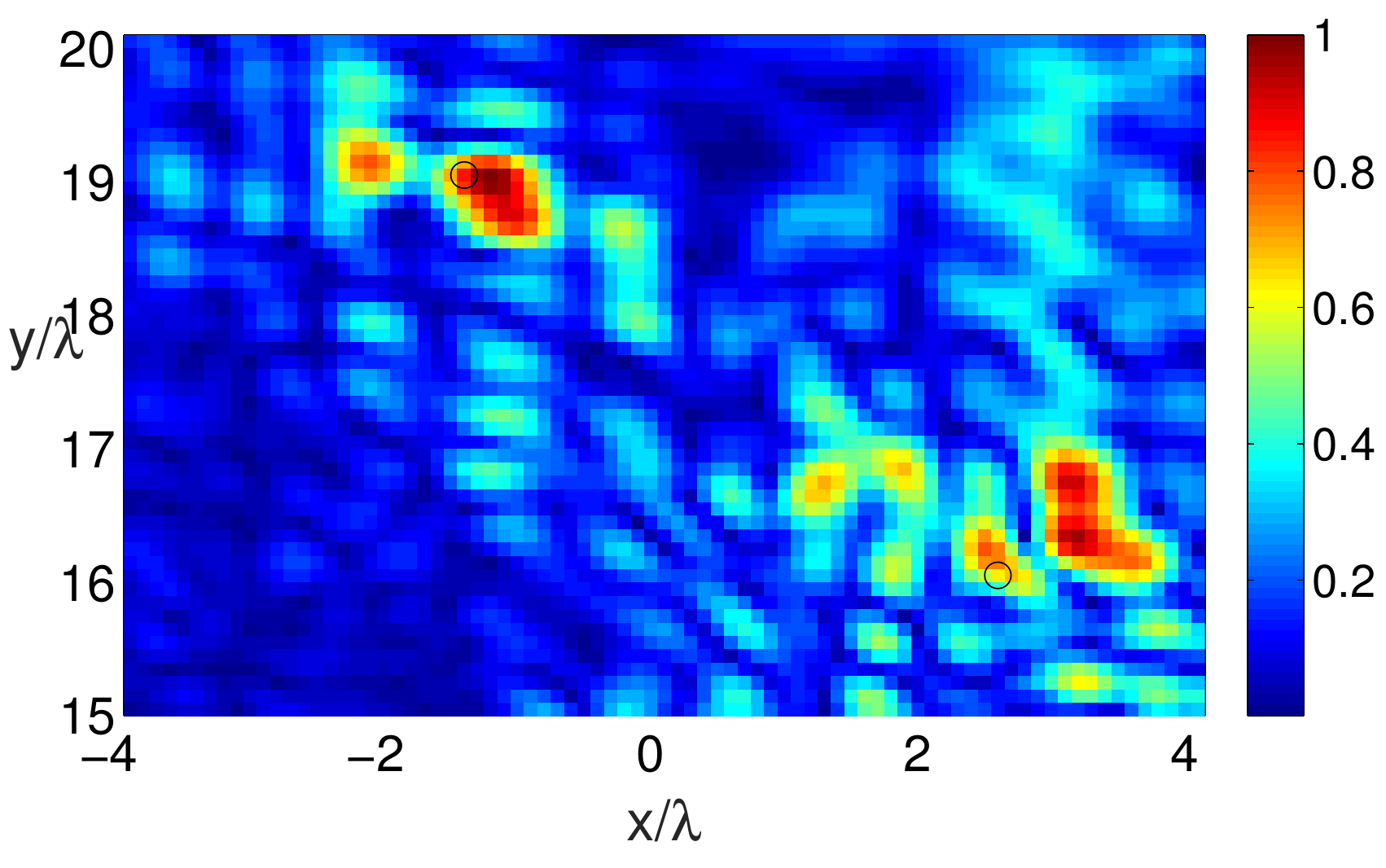}
 \\
 \vspace{-0.in}
\includegraphics[width=0.4\linewidth]{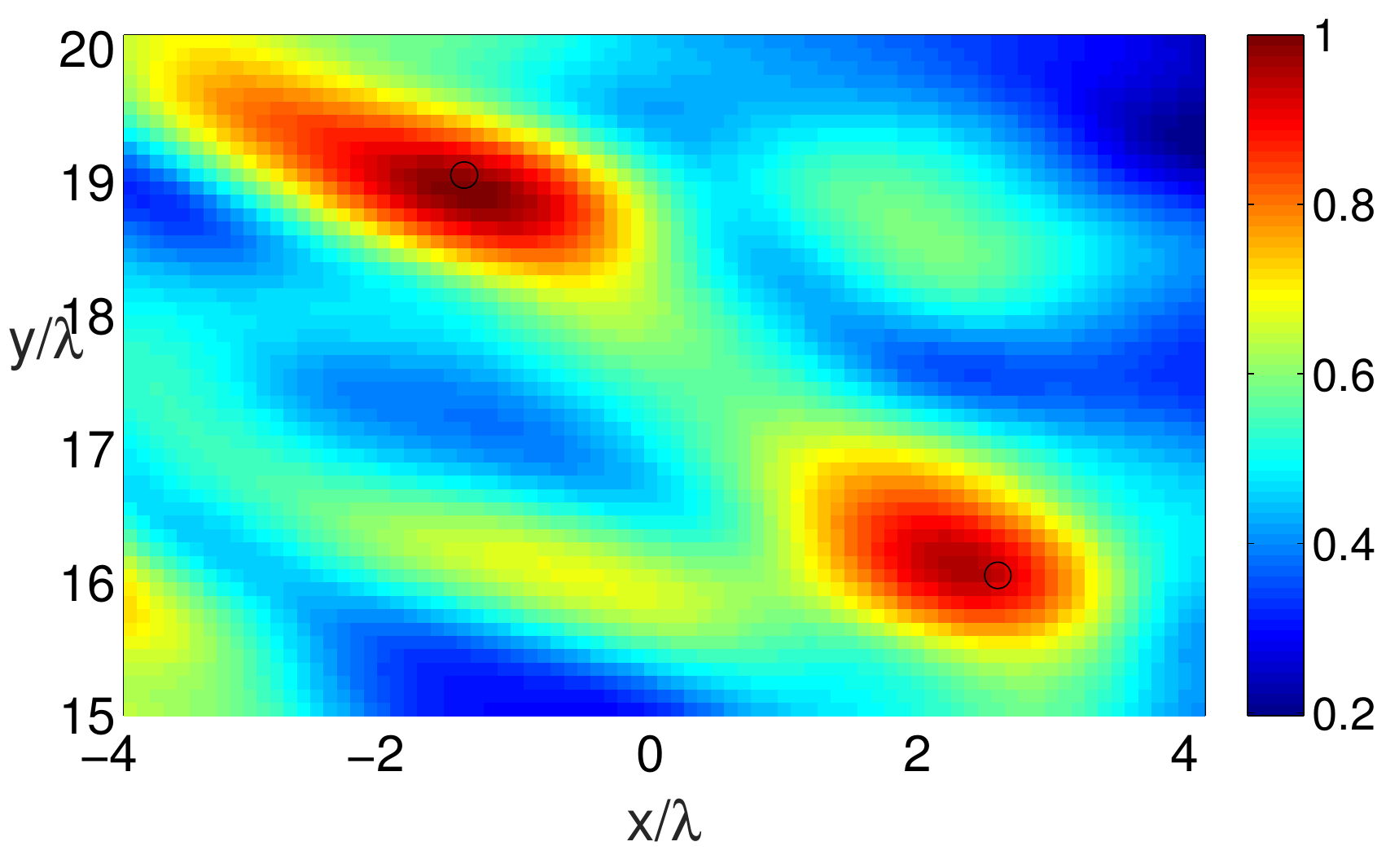}
 \includegraphics[width=0.4\linewidth]{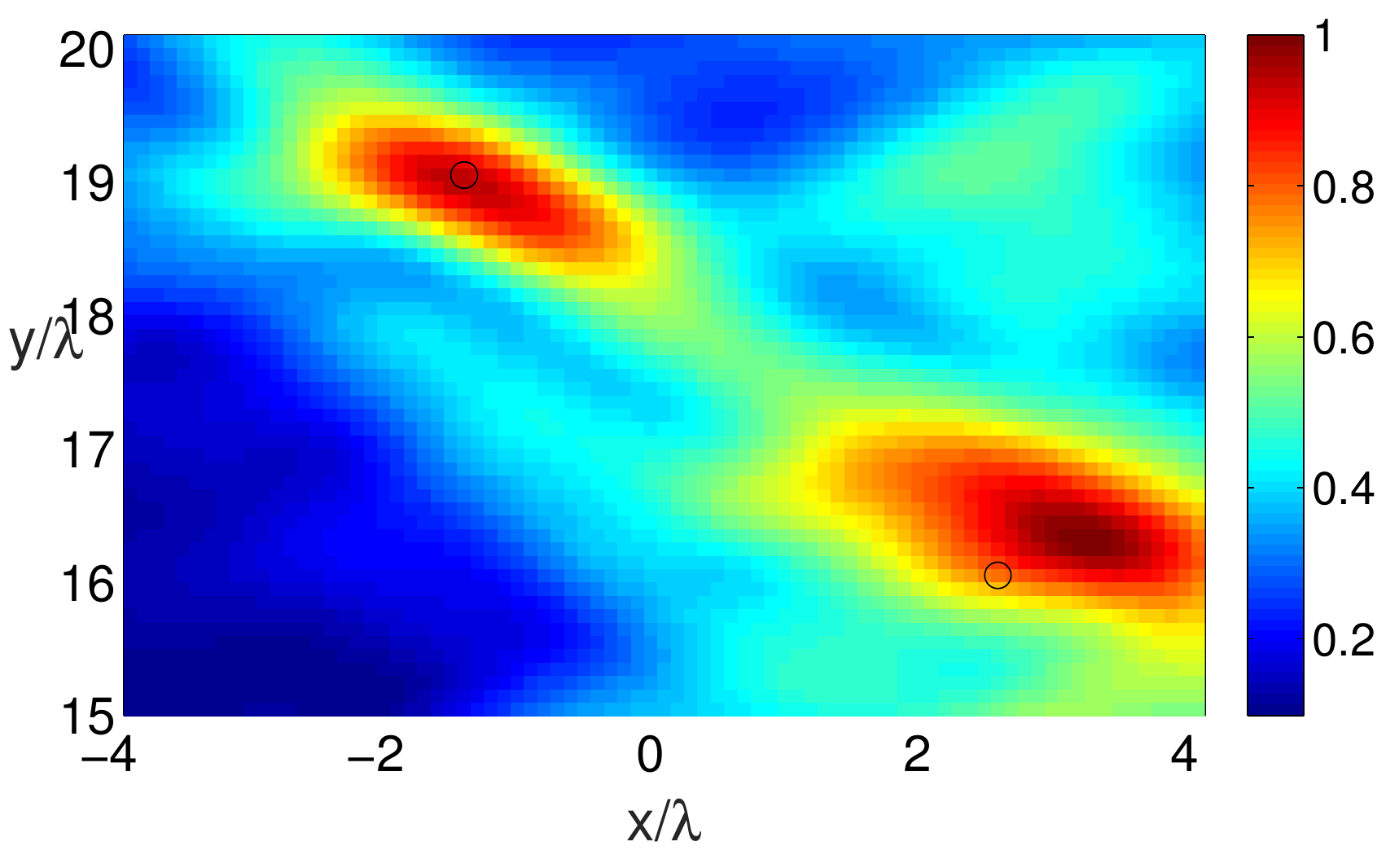}
\vspace{-0.1in}
\caption{The migration (top row) and CINT (bottom) row images of the quadratic susceptibility $\eta_2$ in two realizations of the random medium.  The scatterer location is shown by a black circle.}
\label{MC2}
\end{figure}

The images of the quadratic susceptibility, displayed in Figure \ref{MC2} are similar to those of the linear susceptibility, 
shown in Figure \ref{MC2}. This is because we limited the search domain to the vicinity of the scatterer, and as predicted by the theory
in section \ref{sect:anal}, the uncompensated incident waves in the array data do not have an effect far from the array. However, as shown in Figure \ref{entire1}, these waves cause strong artifacts of the images of the linear susceptibility over larger imaging domains. The image of the quadratic susceptibility is not affected by the direct waves and is clearly better.

\begin{figure}[t]
\centering
\includegraphics[width=0.4\linewidth]{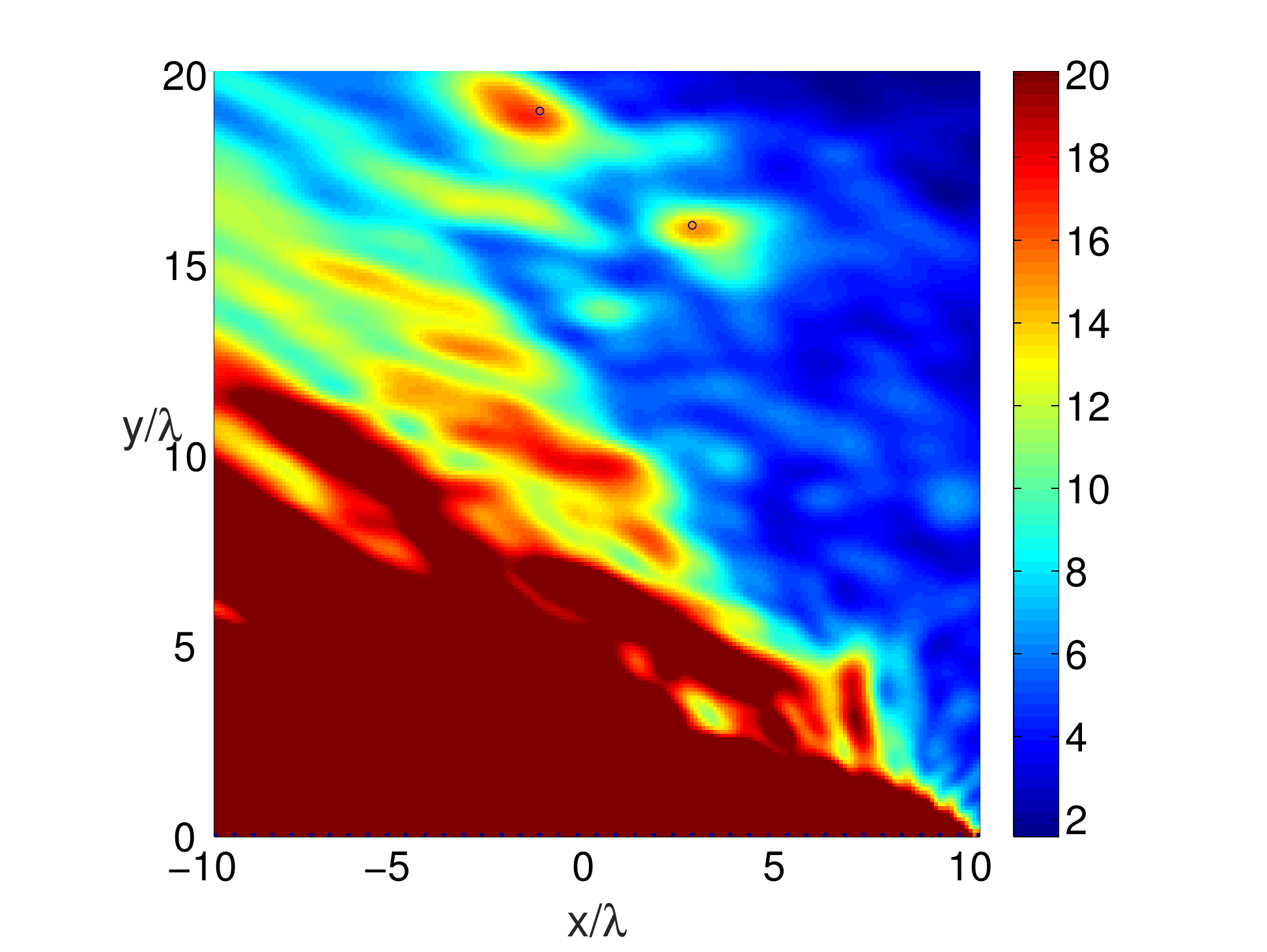}
\includegraphics[width=0.4\linewidth]{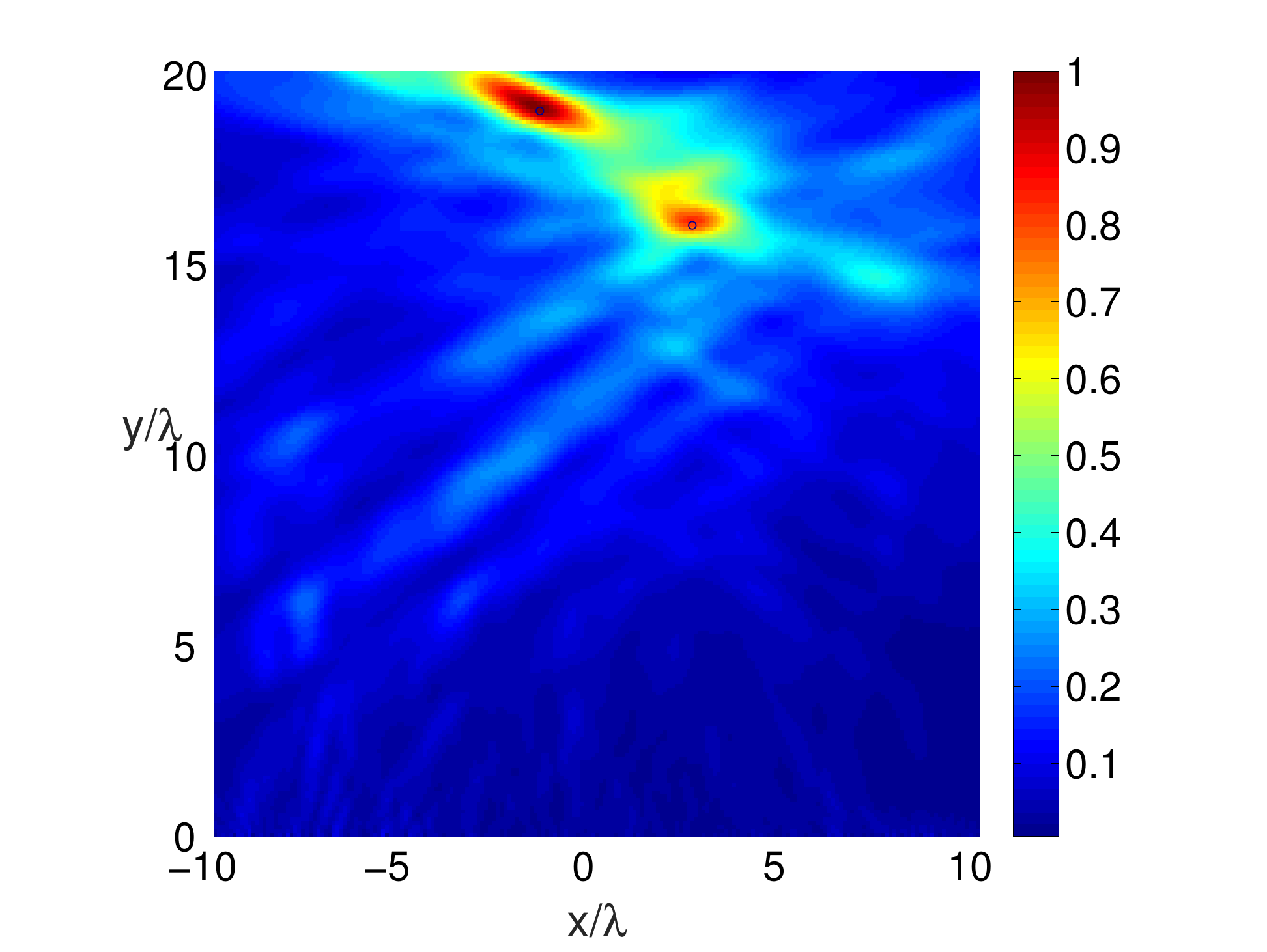}
\vspace{-0.1in}
\caption{CINT images of the linear susceptibility (left) and quadratic susceptibility (right). Here we used 81 incident directions and 81 detectors.}
\label{entire1}
\end{figure}

\section{Conclusions}

We have investigated the problem of optical imaging of small scatterers with second-harmonic light in random media.
We have found that the performance of CINT is superior to that of migration imaging. That is, images obtained by CINT are more robust to statistical fluctuations in the background medium. This observation is consistent with a resolution analysis that is carried out using a geometrical optics model for light propagation in random media. It is also in accord with numerical simulations for which the assumptions of the geometrical optics model do not hold.

\section*{Acknowledgements}

L.B. was partially supported by the NSF grant DMS-1510429. A.M. was partially supported by the NSF grant DMS-1619821 and the University of Houston New 
Faculty Research Program. J.C.S. was partially supported by the NSF grants DMS-1619907 and DMR-1120923. 

\appendix
\section{Statistical moments}
\setcounter{equation}{0}
In this appendix we calculate the statistical moments of the wave
fields.  We begin in section \ref{app:secondmom} with the second
moments of the random phases, which are approximately Gaussian
distributed in our scaling regime described in section \ref{sect:anal.0}. Then we calculate the second moments 
of the wave fields in section \ref{ap:lem2}.

\subsection{Moments of the random phases}
\label{app:secondmom} 
We establish first the following result:
Let $\bx, \bx'$ be two points in the set  $\mathbb{A}$ defined in \eqref{eq:ap3D}. The second moments of
the  processes \eqref{eq:A6} and \eqref{eq:A8} are approximated by
\begin{align}
  \mathbb E[\nu(\bx,\by)\nu(\bx',\by)] \approx \frac{\sqrt{2 \pi}
    \sigma^2 \ell |\bx'-\by|}{4} \int_0^1 dt\, e^{- \frac{t^2}{2
      \ell^2} \left| \bx'_\perp-\bx_\perp \right|^2},
    \label{eq:A17}
\end{align}
and 
\begin{align}
  \mathbb E[\gamma(\bx,\bth)\gamma(\bx',\bth')] \approx \frac{\sqrt{2
      \pi} \sigma^2 \ell |\bx'-\bx'^{(i)}(\bth')|}{4} \int_0^1 dt \,
  e^{- \frac{1}{2 \ell^2} \left|\bm{P}_{\bvth}\left[
      (1-t)(\bx'^{(i)}(\bth')-\bx^{(i)}(\bth)) + t
      (\bx'-\bx)\right]\right|^2}.
    \label{eq:A17g}
\end{align}
The cross-moments satisfy
\begin{equation}
\left| \mathbb E[\gamma(\bx,\bth) \nu(\bx',\by)] \right| \ll
\la^2, \quad \left| \mathbb E[\gamma(\by,\bth)
  \nu(\bx',\by)] \right| \ll \la^2.
\label{eq:A18}
\end{equation}

\vspace{0.05in}
To derive equation \eqref{eq:A17}, we use \eqref{eq:A6} and the Gaussian
autocorrelation \eqref{eq:A2} to obtain
\begin{align}
\mathbb E[\nu(\bx,\by)\nu(\bx',\by)]
&=\frac{\sigma^2|\bx-\by||\bx'-\by|}{4}\int_0^1dt\int_0^1dt'\,
\mathbb E\left[\mu
  \left(\frac{(1-t)\by}{\ell}+\frac{t\bx}{\ell}\right)
  \mu\left(\frac{(1-t')\by}{\ell}+
  \frac{t'\bx'}{\ell}\right)\right] \nonumber
\\ &=\frac{\sigma^2|\bx-\by||\bx'-\by|}{4}\int_0^1dt\int_0^1dt'\, e^{
  - \frac{1}{2 \ell^2} \left| (t'-t)(\bx-\by) +  t'
    (\bx'-\bx)\right|^2 },
\label{eq:Ap1}
\end{align}
for $\bx,\bx' \in A$. We change variables
\[
(t,t') \to (\tilde{t},t'), \quad \tilde t = (t'-t) |\bx-\by|/\ell,
\]
with $t' \in (0,1)$ and $\tilde{t} \in
(-(1-t')|\bx-\by|/\ell,t'|\bx-\by|/\ell)$, and use that 
$
|\bx-\by|/\ell \approx L/\ell \gg 1
$
 to extend the $\tilde{t}$ integral to the real line.
We obtain
\begin{align}
  \mathbb E[\nu(\bx,\by)\nu(\bx',\by)] \approx \frac{\sqrt{2 \pi}
    \ell |\bx'-\by|}{4} \int_0^1 dt' e^{- \frac{1}{2 \ell^2}
    \left|t'\bm{P} (\bx'-\bx)\right|^2},
    \label{eq:Ap2}
\end{align}
where $\bm{P}$ is the orthogonal projection on the plane orthogonal to
$\bx-\by$.  In our case we have
\begin{equation}
  \left| \frac{\bx-\by}{|\bx-\by|}-\bn \right| =
  O\left(\frac{a}{L}\right) \ll 1,
  \label{eq:Ap3}
\end{equation}
and we can estimate the
projection in \eqref{eq:Ap2} by 
\begin{equation}
   \frac{\bm{P}(\bx'-\bx)}{\ell} = \frac{\bm{P}_{_{\bn}}(\bx'-\bx)}{\ell}
   + O\left(\frac{a^2}{\ell L}\right), \quad \bm{P}_n = I - \bn \bn^T.
   \label{eq:Ap4}
\end{equation}
The residual in \eqref{eq:A4} is negligible by assumption
\eqref{eq:as1}, that gives
\[
\frac{a^2}{\ell L} \ll \frac{\sqrt{\la L}}{\ell} \ll 1.
\]
Equation \eqref{eq:A17} follows from \eqref{eq:Ap2} and \eqref{eq:Ap4}.

\vspace{0.05in}
The derivation of \eqref{eq:A17g} is essentially the same, so let us
calculate the cross-moments. We obtain from \eqref{eq:A6},
\eqref{eq:A8} and the Gaussian autocorrelation \eqref{eq:A2} that
\begin{align}
\mathbb E[\gamma(\bx,\bth)\nu(\bx',\by)]
&=\frac{\sigma^2|\bx-\bx^{(i)}(\bth)||\bx'-\by|}{4}\int_0^1dt\int_0^1dt'\,
e^{-\frac{1}{2 } \big|\frac{t'(\bx-\bx^{(i)}(\bth))}{\ell} +
  \frac{t(\bx'-\by)}{\ell} + \frac{\by-\bx}{\ell}\big|^2
},
\label{eq:Ap5}
\end{align}
where the result is obviously positive, so no absolute value is
needed.  Changing variables
\[
(t,t') \to (T,T'), \quad T = t \frac{|\bx'-\by|}{\ell}, \quad T'
= t' \frac{|\bx-\bx^{(i)}(\bth)|}{\ell},
\]
and extending the integrals to the real line, we obtain the upper bound
\begin{align}
\mathbb E[\gamma(\bx,\bth)\nu(\bx',\by)] &\le
\frac{\sigma^2\ell^2}{4}\int_{-\infty}^\infty dT \int_{-\infty}^\infty
dT'\, e^{-\frac{1}{2 } \left|T' \bth + T
  \frac{(\bx'-\by)}{|\bx'-\by|} + \frac{\by-\bx}{\ell}\right|^2
}.
\label{eq:Ap6o}
\end{align}
Expanding the square in the exponent,
\begin{align}
\left|T' \bth + T \left(\frac{\bx'-\by}{|\bx'-\by|}\right)
+ \frac{\by-\bx}{\ell}\right|^2 &= \left[ T' + T
  \frac{\bth \cdot (\bx'-\by)}{|\bx'-\by|}+
  \frac{\bth\cdot (\by-\bx)}{\ell}\right]^2 \notag\\
  &+
\left| \bm{P}_{\bth} \left[T
  \frac{(\bx'-\by)}{|\bx'-\by|} +
  \frac{\by-\bx}{\ell}\right]\right|^2\notag
\end{align}
with $\bm{P}_{\bth} = I - \bth \bth^T$, we obtain after integrating in $T'$ that 
\begin{align}
\mathbb E[\gamma(\bx,\bth)\nu(\bx',\by)] &\le \frac{\sqrt{2 \pi}
  \sigma^2\ell^2}{4}\int_{-\infty}^\infty dT e^{-\frac{1}{2} \left|
    \bm{P}_{\bth} \left[T \frac{(\bx'-\by)}{|\bx'-\by|} +
      \frac{\by-\bx}{\ell}\right]\right|^2 }.
\label{eq:Ap6}
\end{align}
To evaluate the $T$ integral we proceed similarly, by decomposing the
vector $\bm{P}_{\bth} (\by-\bx)$ in two parts: one along
the vector $\bm{P}_{\bth}(\bx'-\by)$ and the other orthogonal to it. 
Then, we expand the square in \eqref{eq:Ap6} and obtain after integration in 
$T$ the upper bound 
\begin{align}
\mathbb E[\gamma(\bx,\bth)\nu(\bx',\by)] &\le \frac{2 \pi 
  \sigma^2\ell^2}{4\left|\bm{P}_{\bth}\frac{(\bx'-\by)}{|\bx'-\by|}\right|}.
\label{eq:Ap7}
\end{align} 
Equation \eqref{eq:A18} follows from the assumptions
\eqref{eq:APA2}, $\bvth \perp \bn$  and equation  \eqref{eq:Ap4} which give
\[
\left|\bm{P}_{\bth}\frac{(\bx'-\by)}{|\bx'-\by|}\right| = O(1),
\]
so that the right hand side in \eqref{eq:Ap7} is $O(\sigma^2 \ell^2)$. 
But by assumption \eqref{eq:as2},
\[
\sigma^2 \ell^2 \ll \frac{\la \ell^3}{L^2} \ll {\la^2},
\]
with the last inequality implied by the upper bound on $\ell$ in \eqref{eq:as1}.

\subsection{Second moments of the wave fields}
\label{ap:lem2}
Let us suppose, without loss of generality, that 
\begin{equation}
|\bx-\by| \ge |\bx'-\by|, \quad |\bx-\bx^{(i)}(\bth)| \ge
|\bx'-\bx'^{(i)}(\bth')|,
\label{eq:AS1}
\end{equation}
and use the results of the preceding section to derive the second moments
of the wave fields. 

\subsubsection{Derivation of moment formula \eqref{eq:MOMG}:} Using that the phases are approximately Gaussian,
we obtain 
\begin{align}
\EE\left[ G(\bx,\by;j\om){G^\star(\bx',\by;j \om)}\right] \approx
G_0(\bx,\by;j\om){G_0^\star(\bx',\by;j \om)}\nonumber \\
\times  \exp\Big\{-\frac{(j
  k)^2}{2} \EE\Big[\big(\nu(\bx,\by)-\nu(\bx',\by)\big)^2\Big]\Big\},
\label{eq:A19}
\end{align}
with exponent  written using \eqref{eq:A17} and definition 
\eqref{eq:A15} as follows
\begin{align*}
\frac{(j k)^2}{2}
\EE\big\{\big[\nu(\bx,\by)-\nu(\bx',\by)\big]^2\big\} =
\frac{|\bx-\by|}{\ell_j^s} + \frac{|\bx'-\by|}{\ell_j^s} - \frac{2
  |\bx'-\by|}{\ell_j^s} \int_0^1 dt \,\exp \Big(- \frac{t^2}{2 \ell^2} \left|\bx'_\perp-\bx_\perp\right|^2\Big).
\end{align*}
Since $|\bx-\by|, |\bx'-\by| \approx L \gg \ell_j^s$, this is very large and
therefore \eqref{eq:A19} is negligible, unless the integral is close
to one. This happens when
\begin{equation}
\big|\bx'_\perp-\bx_\perp\big| \ll \ell,
\label{eq:SMAS}
\end{equation}
in which case we can use a Taylor expansion of the exponential and
obtain the approximation
\begin{align*}
\int_0^1 dt \,e^{- \frac{t^2}{2 \ell^2} \left| \bx'_\perp-\bx_\perp\right|^2} \approx 1 -\frac{\big|\bx'_\perp-\bx_\perp\big|^2 }{6 \ell^2}.
\end{align*}
Substituting in \eqref{eq:A19}, we get 
\begin{align}
\EE\left[ G(\bx,\by;j\om){G^\star(\bx',\by;j \om)}\right] \approx
G_0(\bx,\by;j\om){G_0^\star(\bx',\by;j \om)} \\
\times \exp\Big(-\frac{\left||\bx-\by|-|\bx'-\by|\right|}{\ell_j^s}-\frac{|\bx'_\perp-\bx_\perp|^2 }{2 X_{d,j}^2}\Big),
\label{eq:A19p}
\end{align}
with 
\begin{equation}
X_{d,j} = \ell \sqrt{\frac{3 \ell_j^s}{2 |\bx'-\by|}} \approx \ell
\sqrt{\frac{3 \ell_j^s}{2 L}}.
\label{eq:A19def}
\end{equation}
We also have 
\begin{equation}
|\bx-\by| - |\bx'-\by| \approx
\frac{(\bx-\by)}{|\bx-\by|} \cdot (\bx-\bx') = \bn \cdot
(\bx - \bx') + O\left(\frac{a^2}{L}\right) =
O\left(\frac{a^2}{L}\right),
\label{eq:A20}
\end{equation}
where we used \eqref{eq:Ap3} and that the array aperture  is orthogonal to
$\bn$. But definition \eqref{eq:A15} of $\ell_j^s$ and assumptions
\eqref{eq:as1}-\eqref{eq:as2} give 
\begin{equation}
\frac{\left||\bx-\by| - |\bx'-\by|\right|}{\ell_j^s} =
O\left(\frac{a^2}{L \ell_j^s}\right) = O \left(\frac{a^2\ell \sigma^2}{L
  \la^2}\right) \ll \frac{a^2 \ell^2}{\la L^3} \ll \frac{a^4}{\la L^3}
\ll 1, \label{eq:estimaterange}
\end{equation}
so the first term in the exponent of \eqref{eq:A19p} is negligible.
Equation \eqref{eq:MOMG} follows from \eqref{eq:A19p}, \eqref{eq:A20}
and the paraxial approximation \eqref{eq:APAR3} of $G_0$. This formula
is derived under assumption \eqref{eq:SMAS}. If this doesn't hold, the
moment is exponentially small, of the order $\exp(-2 L/\ell_j^s)$, 
as explained above. This is captured in
the expression \eqref{eq:A20} by the exponential decay on the scale
$X_{d,j}$, which is much smaller than $\ell$ because $\ell_j^s \ll
L$. 

\subsubsection{Derivation of moment formula \eqref{eq:MOMDi}:} 
Again, using the
approximate Gaussian nature of the phases, we obtain from definition
\eqref{eq:A7} that 
\begin{align}
\EE\left[u_1^{(i)}(\bx,\bth){u_1^{(i)}(\bx',\bth')}^*\right] \approx
\exp \Big\{i k (\bx \cdot \bth - \bx' \cdot \bth') - \frac{k^2}{2}
  \EE\Big[\left(\gamma(\bx,\bth)-\gamma(\bx',\bth')\right)^2\Big]\Big\},
\label{eq:ApA1}
\end{align}
with the last term in the exponent following from \eqref{eq:A17g}
\begin{align}
\frac{k^2}{2}\EE\left\{\left[\gamma(\bx,\bth)-\gamma(\bx',\bth')\right]^2\right\}
\approx \frac{|\bx-\bx^{(i)}(\bth)|}{\ell_j^s} +
\frac{|\bx'-\bx'^{(i)}(\bth')|}{\ell_j^s} - \frac{2
  |\bx'-\bx'^{(i)}(\bth')|}{\ell_j^s} \nonumber \\ \times \int_0^1 dt \,
  \exp\Big[- \frac{1}{2 \ell^2} \left|\bm{P}_{\bvth}\left[
      (1-t)(\bx'^{(i)}(\bth')-\bx^{(i)}(\bth)) + t
      (\bx'-\bx)\right]\right|^2\Big].\label{eq:ApA1p}
\end{align}
We conclude as above that since $|\bx-\bx^{(i)}(\bth)|,
|\bx'-\bx'^{(i)}(\bth)| \gg \ell_j^s$, the right hand side in \eqref{eq:ApA1}
is small unless 
\begin{equation}
\frac{|\bm{P}_{\bvth}(\bx'-\bx)|}{\ell} \ll 1, \quad
\frac{|\bm{P}_{\bvth}(\bx'^{(i)}(\bth')-\bx^{(i)}(\bth))|}{\ell} \ll 1.
\label{eq:ApA2}
\end{equation}
By definition 
\begin{equation}
\bx^{(i)}(\bth) = \bx - |\bx-\bx^{(i)}(\bth)| \bth, \quad
\bx'^{(i)}(\bth') = \bx' - |\bx'-\bx'^{(i)}(\bth')| \bth', 
\label{eq:ApA3}
\end{equation}
so the last inequality in \eqref{eq:ApA2} implies that the angle between $\bth$ and
$\bth'$ must be small. 

With the assumption \eqref{eq:ApA2}, we can approximate 
the integral using the Taylor expansion of the exponential,
\begin{align}
\exp\Big\{- \frac{k^2}{2}
  \EE\Big[\Big(\gamma(\bx,\bth)-\gamma(\bx',\bth')\Big)^2\Big]\Big\}
\approx \exp\left[-\frac{\left||\bx-\bx^{(i)}(\bth)| -
    |\bx'-\bx'^{(i)}(\bth')|\right|}{\ell_1^s} \right. \nonumber \\
    \left. - \frac{|\bm{P}_{\bvth} \tbx|^2
    + \tbx \cdot \bm{P}_{\bvth} \tbx^{(i)} + |\bm{P}_{\bvth}
    \tbx^{(i)}|^2}{2 X_{d,1}^2}\right],
\label{eq:ApA4}
\end{align}
with $\ell_1^s$ defined in \eqref{eq:A15}, $X_{d,1}$ defined in \eqref{eq:A19def}, 
$
\tbx = \bx-\bx',$ and  $ \tbx^{(i)} = \bx^{(i)}(\bth)-\bx'^{(i)}(\bth').$
The first term in the exponential in \eqref{eq:ApA4} is of the same order as 
in the estimate \eqref{eq:estimaterange}, and is negligible. 
The second term can be rewritten using \eqref{eq:ApA3}, which gives 
\begin{equation}
\tbx^{(i)} = \tbx -
\left[|\bx-\bx^{(i)}(\bth)|-|\bx'-\bx'^{(i)}(\bth')|\right]\oth -
\Big[|\bx-\bx^{(i)}(\bth)|+|\bx'-\bx'^{(i)}(\bth')|\Big]\frac{\tth}{2},
\label{eq:ApA5}
\end{equation}
where  $
\oth = ({\bth + \bth'})/{2} $ and $ \tth = \bth-\bth'.
$
Note that 
\[
\frac{\left||\bx-\bx^{(i)}(\bth)|-|\bx'-\bx'^{(i)}(\bth')|\right|}{X_{d,1}}
= O\left(\frac{a^2}{L X_{d,1}}\right) = O \left(\frac{a^2 \sigma}{\la
  \sqrt{\ell L}}\right) \ll \frac{a^2}{\sqrt{\la L^3}} \ll 1,
\]
where we used the definition of $X_{d,1}$ and the scaling assumptions
\eqref{eq:as1}-\eqref{eq:as2}. Thus, we 
can neglect the second term in the right hand-side of \eqref{eq:ApA5}, 
and get 
\begin{equation}
\tbx^{(i)} \approx \tbx - |\bx-\bx^{(i)}(\bth)| \tth.
\label{eq:ApA5n}
\end{equation}
Gathering the results and substituting in \eqref{eq:ApA1} and \eqref{eq:ApA4} we obtain 
\eqref{eq:MOMDi}. 

\subsubsection{Derivation of moment formula \eqref{eq:Mig2}:} 
This follows from 
\begin{align}
\EE\left[ G(\bx,\by;j \om) e^{i j k \gamma(\by,\bth)}
  \right] &\approx G_0(\bx,\by;j \om)\exp\Big\{ -\frac{(j k)^2}{2} \EE
  \left[ \left(\nu(\bx,\by)-\gamma(\by,\bth)\right)^2 \right]\Big\}
\nonumber \\ &\approx G_0(\bx,\by;j
\om)\exp\left[-\frac{|\bx-\by|}{\ell_j^s} -
  \frac{|\by-\by^{(i)}(\bth)|}{\ell_j^s}\right],
\end{align}
where the first approximation is because the phases are approximately
Gaussian, and the second approximation is because by \eqref{eq:A18},
\begin{equation}
k^2 \left|\EE \left[\nu(\bx,\by)\gamma(\by,\bth)\right]\right| \ll 1.
\label{eq:factoriz}
\end{equation}

\section{Statistics of the migration image}
\label{ap:mig}
\setcounter{equation}{0}
To find the expectation $\EE\big[\IMj(\by)\big]$ of the migration imaging function 
and evaluate the SNR, we use the following moment factorizations implied by \eqref{eq:factoriz},
\begin{align}
\EE\left[ G(\bx,\by;j \om) {G^\star(\bx',\by;j\om)}
  e^{i j k [\gamma(\by,\bth)-\gamma(\by,\bth')]} \right] \approx \EE\left[ G(\bx,\by;j \om)
  {G^\star(\bx',\by;j\om)}\right]  \nonumber \\
  \times  \EE
\left[e^{i j k [\gamma(\by,\bth)-\gamma(\by,\bth')]}
  \right], \label{eq:Mig3}
\end{align}
and 
\begin{align}
\EE\left[{G^\star(\bx',\by;j\om)}
  e^{i j k [\gamma(\by,\bth')-\gamma(\bx,\bth)]} \right]
\approx \EE\left[ 
  {G^\star(\bx',\by;j\om)}\right]\EE
\left[e^{i j k [\gamma(\by,\bth')-\gamma(\bx,\bth)]}
  \right]. \label{eq:Mig4}
\end{align}
We analyze separately the imaging of the linear and quadratic susceptibilities.

\subsection{Imaging of the linear susceptibility}
Definitions \eqref{eq:dat1} and \eqref{eq:Mig1}, and the estimates 
 \eqref{eq:A16}, \eqref{eq:Mig2} give that the expectation of the image is 
\begin{align}
&\EE\left[\LM(\bys)\right] \approx - \int_{A} d \bx_\perp \int_{C} d
\bth \, {G_0^\star(\bx,\bys;\om)} \exp\Big[i k \bth \cdot(\bx-\bys)\Big]  \nonumber \\ &+ \int_{A} d \bx_\perp \int_{C} d
\bth \, {G_0^\star(\bx,\bys;\om)} \exp\Big[i k \bth \cdot(\bx-\bys)-\frac{|\bx-\bx^{(i)}(\bth)|}{\ell_1^s}\Big] \nonumber \\&+
k^2 \left< \eta_1 \right> \int_{A} d \bx_\perp \int_{C} d \bth \,
G_0(\bx,\by;\om) {G_0^\star(\bx,\bys;\om)}\exp \Big[i k \bth \cdot
  (\by-\bys)- \frac{(|\bx-\by|+|\by-\by^{(i)}(\bth)|)}{\ell_1^s}\Big].
\label{eq:Mig5}
\end{align}
All the terms but the first in this expression are exponentially small. But
even this term gives a small contribution because of the large phase
\[
k \left[|\bx-\bys| - \bth \cdot (\bx-\bys)\right] = k |\bx-\bys| \left[ 1 -
\bth \cdot \frac{(\bx-\bys)}{|\bx-\bys|} \right] = O (L/\la) \gg 1,
\]
where we used the expression of $G_0(\bx,\bys,\om)$ and that $\bvth \perp \bn$. 

The second moment of the imaging function at the scatterer location is 
\begin{align*}
\EE\left[\left|\LM(\by)\right|^2\right] &\approx \int_{A} d \bx_\perp
  \int_{A} d \bx'_\perp \int_{C} d \bth \int_{C} d
  \bth' \, G_0(\bx',\by;\om) {G_0^\star(\bx,\by;\om)} \nonumber \\
  & 
 \times  \Big\{\exp\Big[i k \Big(\bth \cdot \bx- \bth' \cdot \bx' - \by \cdot
      (\bth-\bth')\Big)\Big] \left[ \EE \left[ \exp \Big(i k
        [\gamma(\bx,\bth)-\gamma(\bx',\bth')]\Big)\right] + 1 \right] \nonumber \\
      &  +
  (k^2\left< \eta_1 \right>)^2 \EE\left[G(\bx,\by;\om)
    {G^\star(\bx',\by;\om)}\right] \EE \left[ \exp\Big(i k
      (\gamma(\by,\bth)-\gamma(\by,\bth'))\Big)\right]\Big\},
\end{align*}
where we dropped all the exponentially small terms. Using (\ref{eq:MOMG}) and
(\ref{eq:MOMDi}), we see that  the result is clearly much larger than the square of 
\eqref{eq:Mig5}. Thus, the standard deviation of the image 
\[
\mbox{std}[\LM(\by)] =
\sqrt{\EE\left[\left|\LM(\by)\right|^2\right]-
  \left|\EE\left[\LM(\by)\right]\right|^2} \approx
\sqrt{\EE\left[\left|\LM(\by)\right|^2\right]}
\]
is much larger than its mean. This gives the small SNR in equation \eqref{eq:expsmall}.

\subsection{Imaging of the quadratic susceptibility}
We obtain similarly from definitions \eqref{eq:dat2} and
\eqref{eq:Mig1}, and the estimates \eqref{eq:A16} and  \eqref{eq:Mig2}
that 
\begin{align}
\EE\left[\QM(\bys)\right] &\approx 4 k^2 \left< \eta_2 \right> \int_{A} d \bx_\perp
\int_{C} d \bth \, G_0(\bx,\by;2\om)
{G_0^\star(\bx,\bys;2\om)}\nonumber \\&
\exp\left[i 2k \bth \cdot (\by-\bys)-
  \frac{(|\bx-\by|+|\by-\by^{(i)}(\bth)|)}{\ell_2^s} \right].
\label{eq:Mig6}
\end{align}
This peaks at $\bys = \by$, where the phase cancells out, but the
peak there is small due to the decaying exponential, because $|\bx-\by|$ and
$|\by-\by^{(i)}(\bth)|$ are much larger than the scattering length $ \ell_2^s$.  The second moment at the scatterer 
location is 
\begin{align}
\EE\left[\left|\QM(\by)\right|^2\right] \approx \big(4 k^2 \left<\eta_2\right>\big)^2\int_{A} d \bx_\perp
\int_{A} d \bx'_\perp \int_{C} d \bth \int_{C} d \bth'
\, G_0(\bx',\by;2\om) {G_0^\star(\bx,\by;2\om)}\nonumber \\
\times \EE\left[G(\bx,\by;2\om)
  {G^\star(\bx',\by;2\om)}\right] \EE \left[ \exp\Big(i 2k
      (\gamma(\by,\bth)-\gamma(\by,\bth'))\Big)\right],
\label{eq:Mig7}
\end{align}
with the expectations in the second line calculated in appendix \ref{ap:lem2}.
These expectations are large for nearby points in the array and nearby directions
of illumination. Substituting in \eqref{eq:Mig7} and comparing
with \eqref{eq:Mig6} leads us to 
\[
\mbox{std}[\QM(\by)] =
\sqrt{\EE\left[\left|\QM(\by)\right|^2\right]-
  \left|\EE\left[\QM(\by)\right]\right|^2} \approx
\sqrt{\EE\left[\left|\QM(\by)\right|^2\right]}.
\]
This gives the small SNR in equation \eqref{eq:expsmall}.

\section{Statistics of the CINT image}
\label{ap:CINT}
\setcounter{equation}{0}
We calculate here the mean and variance of the CINT imaging functions $\ICINTj$, for $j = 1, 2$.
The expression of the mean is needed to quantify the focusing of the image, and the variance is 
needed to assess the robustness  with respect to different realizations of the 
random medium.
\subsection{CINT image of the quadratic susceptibility}
\label{ap:QCint}
The expression of the CINT imaging function is obtained by substituting \eqref{eq:CI1} in \eqref{eq:CI3}
\begin{align}
\QC(\bys) &= (4 k^2 \left< \eta_2 \right>)^2\times \iint_{\cC} d \bth d \tth \, \exp\Big\{-\frac{|\bm{P}_{\bvth}\tth|^2}{2 \Theta^2}+ 
i 2 k [\gamma(\by,\bth + \tth/2)-\gamma(\by,\bth - \tth/2)]\Big\} \notag\\
&\times \exp\Big\{
i 2 k [(\bth + \tth/2) \cdot (\by-\bys)-(\bth - \tth/2) \cdot (\by-\bys)] 
\Big\}  \nonumber \\
&\times\iint_A d \bx_\perp  d \tbx_\perp \, \exp\Big(-\frac{|\tbx_\perp|^2}{2 X^2}\Big)
G \Big(\bx + \frac{\tbx}{2},\by;2 \om\Big) G^\star\Big(\bx-\frac{\tbx}{2},\by;2\om\Big) \nonumber \\
&\times
G_0^\star \Big(\bx + \frac{\tbx}{2},\bys;2\om\Big) G_0\Big(\bx-\frac{\tbx}{2},\bys;2 \om\Big)  \, .
\label{eq:ACI1}
\end{align}
It is given by the product of the integrals over the direction vectors and the detector coordinates. Because of the statistical 
decorrelation stated in (\ref{eq:Mig2}) (see also the estimate \eqref{eq:factoriz}) we can study separately the statistics of these integrals, denoted by 
\begin{align}
\nonumber
\mathcal{J}_{\cA}(\bys) = &\iint_{\cA} d \bx_\perp  d \tbx_\perp e^{-\frac{|\tbx_\perp|^2}{2 X^2}} 
G \Big(\bx + \frac{\tbx}{2},\by;2 \om\Big) G^\star\Big(\bx-\frac{\tbx}{2},\by;2\om\Big) 
G_0^\star \Big(\bx + \frac{\tbx}{2},\bys;2\om\Big) \\
& \times G_0\Big(\bx-\frac{\tbx}{2},\bys;2 \om\Big),
\label{eq:ACI2}
\end{align}
and 
\begin{align}
\mathcal{J}_{\cC}(\bys) = \iint_{\cC} d \bth d \tth_\perp \, e^{-\frac{|\bm{P}_{\bvth} \tth|^2}{2 \Theta^2}+
i 2 k [(\bth + \tth/2) \cdot (\by-\bys)-(\bth - \tth/2) \cdot (\by-\bys)] + 
i 2 k [\gamma(\by,\bth + \tth/2)-\gamma(\by,\bth - \tth/2)]}.
\label{eq:ACI3}
\end{align}

\subsubsection{Expectation of the imaging function}
The integral $\mathcal{J}_A(\bys) $ models the CINT point spread function for  a source at $\by$, and has been studied 
in \cite{borcea2011enhanced}. Its expectation follows easily from (\ref{eq:MOMG}) and the definition \eqref{eq:defSetAS} of the set $\cA$,
\begin{align}
\EE \left[\mathcal{J}_A(\bys)\right] &\approx \frac{1}{L^4}\int_A d \bx_\perp  \int_{\mathbb{R}^2} d \tbx_\perp \exp\Big[-\frac{|\tbx_\perp|^2}{2 X^2_e} + 
i 2 k \tbx_\perp \cdot 
(\bys_\perp-\by_\perp)\Big] \\&\approx \frac{2 \pi a^2 X_e^2}{L^4} \exp\Big[-\frac{1}{2} \left(\frac{2k X_e |\by_\perp-\bys_\perp|}{L}\right)^2\Big].
\label{eq:ACI4}
\end{align} 
Here  we extended the $\tbx_\perp$ integral to 
the whole plane using that $X_e \sim X_{d,2} \ll a$.  
Similarly, using (\ref{eq:MOMDi}) we get 
\begin{align}
&\EE \left[\mathcal{J}_{C}(\bys)\right] \approx \iint_{\cC} d \bth   d \tth_\perp \, \exp\left[-\frac{|\bm{P}_{\bvth} \tth|^2}{2 \Theta^2_e} + i 2 k \tth
\cdot (\by-\bys)\right] \nonumber \\
& \hspace{0.2in}=  \iint_{\cC} d \bth   d \tth_\perp \, \exp\left[-\frac{\big(\tilde \theta_\xi^2 + \tilde \theta_\zeta^2\big)}{2 \Theta_e^2} + i 2 k \tilde \theta_\zeta 
\bzet \cdot (\by-\bys) + 2 i k \tilde \theta_\xi \big(\bet - \tan \varphi \bvth \big) \cdot (\by-\bys)\right]
\label{eq:ACI4p}
\end{align} 
with $\tilde \theta_\xi$, $\tilde \theta_\zeta$, $\varphi$ parametrizing 
$\bth$ and $\tth$ as in equations \eqref{eq:Th1}--\eqref{eq:Th5}. To write the integral over the set $\cC$, we recall from 
\eqref{eq:Th1} and \eqref{eq:Th4} that 
\[
\bth = \bvth \cos \varphi + \bet(\beta) \sin\varphi + O(\Theta_d^2),
\]
with azimuthal angle $\beta$ parametrizing the vectors $\bet(\beta)$ and $\bzet(\beta)$.  
In the calculation of the Jacobian of the transformation we may neglect the residual in this equation, and obtain 
\[
d \bth = \sin \varphi d \varphi d \beta.
\]
We also see from equation \eqref{eq:Th5} that for any given $\varphi$ and $\beta$ we have, using $\varphi = O(\alpha) = O(a/L) \ll 1$, that  
\[
\partial_{\tilde{\theta_\xi}} \tth  \approx \bet(\beta) \quad \mbox{and} \quad \partial_{\tilde{\theta_\zeta}} \tth = \bzet(\beta),
\]
and since $\bet(\beta)$ and 
$\bzet(\beta)$ are orthonormal, we get  
\[
d \tth = d \tilde \theta_\xi d \tilde \theta_\zeta.
\]
The integrals over $\tilde \theta_\xi$ and $\tilde \theta_\zeta$ may be 
extended to the real line, because the Gaussians are negligible outside $\cC$, and the result is 
\begin{align} 
\EE \left[\mathcal{J}_{C}(\bys)\right] \approx 2 \pi \Theta_e^2 \int_0^\alpha d \varphi \, \sin\varphi \int_0^{2 \pi} d \beta 
\, \exp\Big[-\frac{(2 k \Theta_e)^2}{2} \big[ (\by-\bys) \cdot \bzet(\beta)\big]^2 \Big] \nonumber \\
\times \exp \Big\{- \frac{(2 k \Theta_e)^2}{2} \big[ (\by-\bys) \cdot \bet(\beta) - \tan \varphi (\by-\bys) \cdot \bvth\big]^2 
\Big\}.
\label{eq:ACI4Cor1}
\end{align}
We are interested only in the points $\bys$ for which  $\mathcal{J}_{\cA}(\bys)$ is large, so 
\[
| (\by - \bys) \cdot \bvth| =  O\left(\frac{L}{k X_e}\right).
\]
Since $\Theta_e \approx X_e/L$ and  $\varphi \le \alpha \ll 1$, we have 
\[
k \Theta_e \tan \varphi  |(\by - \bys) \cdot \bvth|\le O(\alpha) \ll 1,
\]
and we can neglect the $\varphi$ dependent term in the exponential in \eqref{eq:ACI4Cor1}.  We also note that 
\[
\big[ (\by-\bys) \cdot \bzet(\beta)\big]^2 + \big[ (\by-\bys) \cdot \bet(\beta)\big]^2 = \big| \bm{P}_{\bvth} (\by-\bys)\big|^2
\]
is independent of $\beta$, so we obtain 
\begin{align} 
\EE \left[\mathcal{J}_{C}(\bys)\right] &\approx 2 \pi \Theta_e^2 \exp\Big\{-\frac{\big[ 2 k \Theta_e \big|\bm{P}_{\bvth} (\by-\bys)\big|\big]^2}{2}\Big\}
\int_0^\alpha d \varphi \, \sin\varphi \int_0^{2 \pi} d \beta \nonumber \\
&= 2 \pi^2 \Theta_e^2 \alpha^2 \exp\Big\{-\frac{\big[ 2 k \Theta_e \big|\bm{P}_{\bvth} (\by-\bys)\big|\big]^2}{2}\Big\}.
\end{align}

\subsubsection{Variance of the imaging function}
\label{ap:variance}
The variance of $\mathcal{J}_A$ is calculated  in 
\cite[Appendix E]{borcea2011enhanced}, so we revisit here  the main ideas. The calculation involves the fourth moments of the Green's function 
\eqref{eq:A5}, which are determined by 
\begin{align}
\EE \Big\{ \exp\Big[i 2 k \big[ \nu\big(\bx+\frac{\tbx}{2},\by\big) - \nu\big(\bx-\frac{\tbx}{2},\by\big)-\nu\big(\bx'+\frac{\tbx'}{2},\by\big)+
\nu\big(\bx'-\frac{\tbx'}{2},\by\big)\big]\Big]\Big\} \approx 
e^{- \tau/2},
\label{eq:ACI5}
\end{align}
where we introduced the notation 
\begin{align}
\tau = (2 k)^2 \EE\big\{ \big[ \nu\big(\bx+\frac{\tbx}{2},\by\big) - \nu\big(\bx-\frac{\tbx}{2},\by\big)-\nu\big(\bx'+\frac{\tbx'}{2},\by\big)+
\nu\big(\bx'-\frac{\tbx'}{2},\by\big)\big]^2\big\},
\label{eq:ACI6}
\end{align}
and used the approximate Gaussian distribution of the phases. We need the second moments \eqref{eq:A17}, rewritten as 
\begin{equation}
(2 k)^2 \EE\left[\nu(\bx,\by) \nu(\bx',\by)\right] = \frac{3 \ell^2}{X_{d,2}^2} {h}\left(\frac{|\bx_\perp-\bx'_\perp|}{\ell}\right), \quad 
{h}(z) = \frac{1}{z} \int_0^z dt\, e^{-\frac{t^2}{2}},
\label{eq:ACI7}
\end{equation}
using the definition \eqref{eq:defXd} of the decoherence length. The expression \eqref{eq:ACI6} becomes 
\begin{align*}
\tau = \frac{6 \ell^2}{X_{d,2}^2} \left[ 2 - h\Big(\frac{\tbx_\perp}{\ell}\Big) -  
h\Big(\frac{\tbx'_\perp}{\ell}\Big) + h\Big(\frac{\bx_\perp-\bx'_\perp}{\ell} 
+ \frac{\tbx_\perp + \tbx'_\perp}{2 \ell} \Big) + h\Big(\frac{\bx_\perp-\bx'_\perp}{\ell} 
- \frac{\tbx_\perp + \tbx'_\perp}{2 \ell} \Big)  \right. \\
\left.-h\Big(\frac{\bx_\perp-\bx'_\perp}{\ell} 
+ \frac{\tbx_\perp - \tbx'_\perp}{2 \ell} \Big) 
-h\Big(\frac{\bx_\perp-\bx'_\perp}{\ell} 
- \frac{\tbx_\perp - \tbx'_\perp}{2 \ell} \Big) \right],
\end{align*}
and we can simplify it because $|\tbx_\perp|, |\tbx'_\perp| \lesssim X_{d,2} \ll \ell$, due to the windowing in the 
calculation of the cross-correlations. Expanding in $\tbx_\perp/\ell$ and $\tbx'/\ell'$ we get 
\begin{equation}
\tau \approx  \frac{6}{X_{d,2}^2}\left[\frac{|\tbx_\perp|^2 + |\tbx'_\perp|^2}{6}  + \tbx_\perp
\cdot H \Big(\frac{\bx_\perp-\bx'_\perp}{\ell}\Big) \tbx'_\perp\right],
\label{eq:ACI8}
\end{equation}
where $H$ is the Hessian of $h$, evaluated at $(\bx_\perp-\bx'_\perp)/\ell$.  

Because the Hessian decays,  we note that  the phase
differences at points satisfying $|\bx_\perp-\bx'_\perp| \gg \ell$ are decorrelated
\begin{align*}
\tau \approx \frac{|\tbx_\perp|^2 + |\tbx'_\perp|^2}{X_{d,2}^2} = (2 k)^2 \EE \Big\{\Big[ \nu\Big(\bx+\frac{\tbx}{2},\by\Big) - 
\nu\Big(\bx-\frac{\tbx}{2},\by\Big)\Big]^2 \Big\} +(2 k)^2 \EE \Big\{\Big[ \nu\Big(\bx'+\frac{\tbx'}{2},\by\Big) \\ - 
\nu\Big(\bx'-\frac{\tbx'}{2},\by\Big)\Big]^2 \Big\}. 
\end{align*}
It is only for $|\bx_\perp-\bx'_\perp| \lesssim \ell$ that the Hessian contributes to \eqref{eq:ACI8}. Thus, when calculating 
the variance of the CINT imaging function, we get a contribution only from the set of points 
\[
\{\bx,\bx' \in \mathbb{A}, ~ ~ |\bx_\perp-\bx'_\perp| \lesssim \ell \}.
\]
This is why the SNR of $\mathcal{J}_A$  is of order $a/\ell$. We refer to \cite[Appendix E]{borcea2011enhanced} for more details.

The calculation of the variance of $\mathcal{J}_{C}$ is similar, and the SNR is of the same order.

\subsection{CINT image of the linear susceptibility}
\label{ap:CINTLin}
The expression of the imaging function is obtained by substituting  \eqref{eq:CI6} in \eqref{eq:CI8},
\begin{align}
&\LC(\bys)= \iint_{\cA} d \bx_\perp d \tbx_\perp \iint_{\cC} d \bth  d \tth \, e^{-\frac{|\tbx_\perp|^2}{2 X^2}-\frac{|\bm{P}_{\bvth} \tth|^2}{2 \Theta^2} - i k \tth \cdot \bys}
G_0^\star \Big(\bx+\frac{\tbx}{2},\bys;\om\Big) G_0 \Big(\bx-\frac{\tbx}{2},\bys;\om\Big)  \nonumber  \\ &\times 
\Big\{ \Big[ e^{-i k \gamma \big(\bx - \frac{\tbx}{2},\bth-\frac{\tth}{2}\big)}-1\Big] e^{-i k \big(\bth-\frac{\tth}{2}\big)\cdot\big(\bx - \frac{\tbx}{2} \big) }
+ k^2 \left< \eta_1 \right> G^\star \Big( \bx-\frac{\tbx}{2},\by;\om \Big) e^{-i k \big(\bth -\frac{\tth}{2}\big) \cdot\by - i k \gamma \big(\by,\bth - \frac{\tth}{2}\big) }\Big\} \nonumber \\
&\times   \Big\{\Big[ e^{i k \gamma \big(\bx + \frac{\tbx}{2},\bth+\frac{\tth}{2}\big)}-1\Big] e^{i k \big(\bth+\frac{\tth}{2}\big)\cdot\big(\bx +\frac{\tbx}{2}\big) } +
k^2 \left< \eta_1 \right> G \Big( \bx+\frac{\tbx}{2},\by;\om \Big) e^{i k \big(\bth +\frac{\tth}{2}\big) \cdot \by + i k \gamma \big(\by,\bth + \frac{\tth}{2}\big)
}\Big\} \, .
\label{eq:CILIN}
\end{align}
Using (\ref{eq:MOMG}) and (\ref{eq:MOMDi}), we obtain the expectation 
\begin{align}
\EE&\Big[\LC(\bys) \Big] \approx \iint_{\cA} d \bx_\perp d \tbx_\perp \iint_{\cC} d \bth d \tth \, e^{-\frac{|\tbx_\perp|^2}{2 X^2}-\frac{|\bm{P}_{\bvth} \tth|^2}{2 \Theta^2}} 
G_0^\star \Big(\bx+\frac{\tbx}{2},\bys;\om\Big) G_0 \Big(\bx-\frac{\tbx}{2},\bys;\om\Big) \nonumber  
\\
& \times e^{-i k \tth \cdot \bys} \Big\{ e^{i k \bth \cdot \tbx + i k \tth \cdot \bx} \Big[ \EE \Big[e^{i k \big[\gamma\big( \bx + \frac{\tbx}{2},\bth+\frac{\tth}{2}\big) - 
\gamma\big( \bx + \frac{\tbx}{2},\bth+\frac{\tth}{2}\big)\big]}\Big] + 1 \Big] \nonumber \\
&+ k^4 \left< \eta_1 \right>^2  e^{i k \tth \cdot \by} \EE\Big[G \Big( \bx+\frac{\tbx}{2},\by;\om \Big) G^\star \Big( \bx-\frac{\tbx}{2},\by;\om \Big) \Big] \EE \Big[ e^{i k \big[\gamma\big( \by,\bth+\frac{\tth}{2}\big) - 
\gamma\big( \by,\bth+\frac{\tth}{2}\big)\big]}\Big] \Big\}, \label{eq:ACI9}
\end{align}
where we have neglected the small terms due to decaying exponentials.
We write the right hand side  as the sum of three terms
\begin{equation}
\EE\Big[\LC(\bys) \Big]  \approx \mathcal{T}_1(\bys) + \mathcal{T}_2(\bys) + \mathcal{T}_3(\bys).
\label{eq:ACI10}
\end{equation}     
The first two terms are due to the uncompensated incident wave, and are given by 
\begin{align}
\mathcal{T}_{1}(\bys) = \frac{1}{L^2}\iint_{\cA} d \bx_\perp d \tbx_\perp  \int_{\cC} d \bth{d}\tth \, 
\exp\Big[i k \tbx \cdot \bth -  \frac{i k \tbx_\perp \cdot (\bx_\perp-\bys_\perp)}{L}  + i k \tth \cdot (\bx-\bys) \Big]
\nonumber \\
\times \exp \Big[ -\frac{|\tbx_\perp|^2}{2 X^2}- \frac{|\bm{P}_{\bvth} \tth|^2}{2 \Theta^2}\Big],
\label{eq:ACI11_T1} 
\end{align}
and 
\begin{align}
&\mathcal{T}_{2}(\bys) = \frac{1}{L^2}\iint_{\cA} d \bx_\perp d \tbx_\perp  \int_{\cC} d \bth{d}\tth \, 
\exp \Big[i k \tbx \cdot \bth -  \frac{i k \tbx_\perp \cdot(\bx_\perp-\bys_\perp)}{L} \big) + i k \tth \cdot (\bx-\bys) \Big]
\nonumber \\
&\times  \exp\Big[-\frac{|\tbx_\perp|^2}{2 X^2} - \frac{|\bm{P}_{\bvth} \tth|^2}{2 \Theta^2}-\frac{|\bm{P}_{\bvth} \tth|^2}{2 [X_{d,1}/|\bx-\bx^{(i)}(\bth)|]^2} - \frac{3}{2} \frac{|\bm{P}_{\bvth} \tbx|^2}{X_{d,1}^2} + 
\frac{3}{2} \frac{\bm{P}_{\bvth} \tbx }{X_{d,1}} \cdot \frac{\bm{P}_{\bvth} \tth}{X_{d,1}/|\bx-\bx^{(i)}(\bth)|}\Big].
\label{eq:ACI11_T2}
\end{align}
Here we used the paraxial approximation \eqref{eq:APAR3} and  moment formula \eqref{eq:MOMDi}.
The third term is similar to the expectation of the imaging function for the quadratic susceptibility $\EE[\QC(\bys)]$, so we write it directly,
\begin{align}
\mathcal{T}_3(\bys) = \frac{(2 \pi)^3}{2} \Big(k^2 \left< \eta_1 \right>^2 \alpha \Theta_e \frac{a X_e}{L^2} \Big)^2 
\exp \Big[-\frac{1}{2} \left(\frac{k X_e |\by_\perp-\bys_\perp|}{L}\right)^2 - \frac{1}{2} \left(k \Theta_e 
|\bm{P}_{\bvth}(\by-\bys)|\right)^2\Big] .
\label{eq:ACI12}
\end{align}

To calculate \eqref{eq:ACI11_T1} we recall the definition \eqref{eq:defSetAS} of the set $\cA$ and the parametrization 
of the set $\cC$ described in equations \eqref{eq:Th1}-\eqref{eq:Th5}. We integrate over $\tbx_\perp$ using that $\tbx$ is orthogonal to $\bn$, 
and 
\[
k \tbx\cdot \bth = k |\bth| \tbx \cdot (\cos \varphi\,  \bvth + \sin \varphi \, \bet(\beta)) = k \tbx_\perp (\cos \varphi\,  \bvth_\perp + \sin \varphi \, \bet_\perp(\beta)) + O (k X_{d,1} \Theta_d^2),
\]
with two dimensional vectors $\bvth_\perp$ and $\bet_\perp$ of components of $\bvth$ and $\bet$ in the plane orthogonal to $\bn$. 
The residual is negligible by equations \eqref{eq:defXd}, \eqref{eq:defThed} and assumption \eqref{eq:as1},
\[
k X_{d,1} \Theta_d^2 = O\left(\frac{X_{d,1}^3}{\la L^2}\right) \ll \frac{\ell^3}{\la L^2} \ll 1.
\]
To integrate over $\tth$, more precisely over $\tilde \theta_\xi$ and $\tilde \theta_\zeta$, we use that 
\[
|\bm{P}_{\bvth} \tth|^2 = \tilde \theta_\xi^2 + \tilde \theta_\zeta^2, \quad 
\tth \cdot (\bx-\by) = \tilde \theta_\xi (\bx-\bys) \cdot \big( \bet - \bvth \tan \varphi\big) + \tilde \theta_\zeta  (\bx-\bys) \cdot \bzet.
\]
We obtain that 
\begin{align}
\mathcal{T}_{1}(\bys)  \approx \frac{(2 \pi)^2 X^2 \Theta^2}{L^2} \int_A d \bx_\perp \int_0^\alpha d \varphi \, \sin \varphi \int_0^\beta d \beta 
e^{-\frac{(k X)^2}{2} \big|  \cos \varphi \, \bvth_\perp- \frac{\bx_\perp-\bys_\perp}{L} + \sin \varphi  \bet_\perp(\beta)\big|^2} \nonumber \\
\times e^{-\frac{(k \Theta)^2}{2} \big\{ \big[ (\bx -\bys) \cdot \bzet(\beta)\big]^2 + \big[ (\bx -\bys) \cdot  \bet(\beta)- \tan \varphi (\bx-\bys) \cdot 
\bvth\big]^2\big\}},
\label{eq:cor1}
\end{align}
and note that the result is exponentially small. The first exponential is small because 
\[
\Big|\cos \varphi \, \bvth_\perp- \frac{\bx_\perp-\bys_\perp}{L} + \sin \varphi  \bet_\perp(\beta)\Big| = \cos \varphi + O\left(\frac{a}{L}\right) \approx 1,
\]
and  by definition \eqref{eq:defXd} and assumptions \eqref{eq:as1}-\eqref{eq:as2}, we have
\[
k X = O \left(\frac{X_{d,1}}{\la}\right) = O\left(\frac{\sqrt{\ell}}{\sigma \sqrt{L}}\right) \gg \sqrt{\frac{L}{\la}} \gg 1.
\]
The second exponential is small because 
\[
k \Theta \sqrt{\big[ (\bx -\bys) \cdot \bzet(\beta)\big]^2 + \big[ (\bx -\bys) \cdot  \bet(\beta)- \tan \varphi (\bx-\bys) \cdot 
\bvth\big]^2} = O(k \Theta L) = O(k X_{d,1}) \gg 1.
\]
The calculation of \eqref{eq:ACI11_T2} is similar, slightly more involved, and the result is exponentially small for points 
$\bys$ in the imaging region $R$.  

The calculation of the variance of  $\LC(\bys)$ is very similar to that described in appendix \ref{ap:variance}. It shows that the incident wave 
has a negligible effect at points $\bys \in R$, due to the large deterministic uncompensated phases. The variance is approximately
equal to that of the useful term in the imaging function, which focuses at the scatterer, and the SNR is of order $(a/\ell)^2$, 
as in the case of quadratic susceptibility. 

It remains to verify that the expectation of the imaging function \eqref{eq:CI8} is large at points $\bys$ near the array. We do this by studying the imaging function \eqref{eq:CILIN} at points $\bys$ outside the small search region $R$.
It suffices to consider only the terms that involve the incident waves, because we know from the analysis above that the the waves that interact with the scatterer at $\by$ contribute to the image only in the vicinity of $\by$. We obtain from \eqref{eq:CILIN} that the contribution of the 
incident waves to the expectation of the image is given by the sum of two terms
\begin{align}
\mathcal{T}_1(\bys) &= \iint_{\cA} d \bx_\perp d \tbx_\perp \iint_{\cC} d \bth  d \tth \, \exp\Big(-\frac{|\tbx_\perp|^2}{2 X^2}-\frac{|\bm{P}_{\bvth} \tth|^2}{2 \Theta^2}\Big)
G_0^\star \Big(\bx+\frac{\tbx}{2},\bys;\om\Big)\notag\\
& \times G_0 \Big(\bx-\frac{\tbx}{2},\bys;\om\Big) \exp\Big(i k \tth \cdot (\bx -\bys) + i k \bth \cdot \tbx\Big)\,, 
\label{eq:CILIND1}
\end{align}
and 
\begin{align}
\mathcal{T}_2(\bys) &= \iint_{\cA} d \bx_\perp d \tbx_\perp \iint_{\cC} d \bth  d \tth \, \exp\Big[-\frac{|\tbx_\perp|^2}{2 X^2}-\frac{|\bm{P}_{\bvth} \tth|^2}{2 \Theta^2}\Big]
G_0^\star \Big(\bx+\frac{\tbx}{2},\bys;\om\Big) G_0 \Big(\bx-\frac{\tbx}{2},\bys;\om\Big)  \nonumber \\
&\times \exp\Big[i k \tth \cdot (\bx -\bys) + i k \bth \cdot \tbx-\frac{|\bm{P}_{\bvth} \tth|^2}{2 [X_{d,1}/|\bx-\bx^{(i)}(\bth)|]^2} - \frac{3}{2} \frac{|\bm{P}_{\bvth} \tbx|^2}{X_{d,1}^2}\Big] \\
&\times \exp\Big[
\frac{3}{2} \frac{\bm{P}_{\bvth} \tbx }{X_{d,1}} \cdot \frac{\bm{P}_{\bvth} \tth}{X_{d,1}/|\bx-\bx^{(i)}(\bth)|}\Big]\, .
\label{eq:CILIND2}
\end{align}
The product of the Green's functions in these expressions is approximated by 
\begin{equation}
G_0^\star \Big(\bx+\frac{\tbx}{2},\bys;\om\Big) G_0 \Big(\bx-\frac{\tbx}{2},\bys;\om\Big) \approx \frac{1}{|\bx-\bys|^2} \exp\Big[i k \tbx \cdot 
\frac{(\bx-\by)}{|\bx-\by|}\Big]\, ,
\label{eq:Adit1}
\end{equation}
because 
\begin{equation}
\label{eq:Adit1_1}
k \Big(\Big|\bx + \frac{\tbx}{2} - \bys \Big| - \Big| \bx - \frac{\tbx}{2} - \bys \Big| \Big) = k \tbx \cdot \frac{(\bx-\bys)}{|\bx-\bys|} + O \left(
\frac{|\tbx|^3 (y^R_{\parallel})^2 |\bx_\perp-\bys_\perp|}{\la |\bx-\bys|^5}\right).
\end{equation}
This is assuming $|\bx-\bys| \gg X$, which holds for a fixed $\bys$ at all $\bx$ in the aperture, with the possible exception of a small set,
of radius of order $X$, which makes a negligible contribution to the integrals in \eqref{eq:CILIND1} and \eqref{eq:CILIND2}. 
Under this condition we see that the residual in \eqref{eq:Adit1_1} is negligible for search points near the array (with small enough $y^R_{\parallel}$), 
and we can approximate the integrals \eqref{eq:CILIND1} and \eqref{eq:CILIND2} using the approximation \eqref{eq:Adit1}.

Substituting \eqref{eq:Adit1} in 
\eqref{eq:CILIND1}, and integrating over $\tbx$ and $\tth$ we get that 
\begin{align} 
\mathcal{T}_1(\bys)  &\approx (2 \pi X \Theta/2)^2 \int_A d \bx_\perp \frac{1}{|\bx-\bys|^2} \int_0^{\alpha} d \varphi \, \sin \varphi \int_0^{2 \pi} d \beta \, 
\exp \Big[-\frac{(k X)^2}{2} \left|\bm{P}_{\bn} \left(\bth - \frac{\bys-\bx}{|\bys-\bx|}\right)\right|^2\Big] \nonumber \\
&\times \exp \Big\{ - \frac{(k \Theta)^2}{2} 
\left\{[(\bx-\bys)\cdot \bzet(\beta)]^2 + [ (\bx-\bys) \cdot (\bet(\beta) - \tan \varphi \bvth)]^2\right\}\Big\},
\label{eq:Adit2}
\end{align}
with $\bth$ parametrized as in equations \eqref{eq:Th1} and \eqref{eq:Th4}. It is difficult to evaluate these integrals explicitly, unless we make 
further scaling assumptions on the location of $\bys$. However, it is clear that \eqref{eq:Adit2} is large when 
\begin{equation}
\left|\bm{P}_{\bn} \left(\bth - \frac{\bys-\bx}{|\bys-\bx|}\right)\right| \lesssim  \frac{1}{k X}  = O\left(\frac{\la}{X_{d,1}}\right) \ll 1,
\label{eq:Adit4}
\end{equation}
for most directions $\bth$ in the cone of illuminations. Equations \eqref{eq:Th1} and \eqref{eq:Th4}, and the assumed orthogonality 
of $\bvth$ and $\bn$  give that 
\[
\bm{P}_{\bn} \bth = \cos \varphi \bvth + \sin \varphi \bm{P}_{\bn} \bet(\beta),
\]
and since $\varphi \le \alpha = O(a/L) \ll 1$, we see that the image is large when 
\begin{equation}
\label{eq:Adit3}
\bvth \cdot \frac{(\bys-\bx)}{|\bys-\bx|} \approx \cos \varphi \approx 1.
\end{equation}
This can hold only at points $\bys$ near the array. The second exponential in \eqref{eq:Adit2} is large when 
\[
\left| \bm{P}_{\bvth} (\bx-\bys) \right| \lesssim \frac{1}{k \Theta} = O\left(\frac{\la L}{X_{d,1}}\right),
\]
which is consistent with \eqref{eq:Adit3}. 

The calculation of \eqref{eq:CILIND2} is similar, and leads to the same conclusion. We end with the remark that 
the set of points where $\mathcal{T}_1(\bys)$ and $\mathcal{T}_2(\bys)$ are large   depends on the aperture size and
the opening angle of the cone of illuminations. Indeed, equation  \eqref{eq:Adit3} gives that the image is large when 
\[
\bvth \cdot \frac{(\bys-\bx)}{|\bys-\bx|} = O(\cos \alpha), 
\]
so the larger $\alpha$ is, the further from the array $\bys$ can be.  Moreover, the larger the aperture is, the more points $\bys$ 
satisfy this equation, for at least some subset of the detector locations in the array. 

\section{Numerical solution of the forward problem}
\label{app:numerics}

To solve the system of nonlinear Helmholtz equations
\eqref{eq:waveu}-\eqref{eq:wavev} numerically we employ the fixed
point iteration described below. We denote by $H_k$ and $H_{2k}$,
respectively the linear operators in
\eqref{eq:waveu}-\eqref{eq:wavev}:
\begin{align}
H_k & = \Delta + k^2 ( 1 + 4\pi\eta(\bx) + 4\pi\eta_1(\bx) ),
\label{eqn:hk} \\
H_{2k} & = \Delta + (2k)^2 ( 1 + 4\pi\eta(\bx) + 4\pi\eta_1(\bx)
).
\label{eqn:h2k}
\end{align}
We also introduce the successive approximations to the solutions of
\eqref{eq:waveu}-\eqref{eq:wavev} as
\begin{align}
u^{(j)}(\bx) & \approx u_1(\bx) - u^{(i)}_1(\bx), \label{eqn:uj} \\
v^{(j)}(\bx) & \approx u_2(\bx). \label{eqn:vj}
\end{align}
We substitute (\ref{eqn:uj})--(\ref{eqn:vj}) into
\eqref{eq:waveu}-\eqref{eq:wavev}, and obtain for $j=0,1,\ldots$ the
following fixed point iteration
\begin{align}
u^{(j+1)} & = - 4 \pi k^2 H_k^{-1} \left[ 
2 \eta_2(\bx) v^{(j)}(\bx) (u^{(j)}(\bx) + u_i(\bx))^*
+ (\eta_0(\bx) + \eta_1(\bx))u_i(\bx)
\right], \label{eqn:fpu} \\
v^{(j+1)} & = - 16 \pi k^2 H_{2k}^{-1} \left[ 
\eta_2(\bx)(u^{(j+1)}(\bx) + u_i(\bx))^2
\right], 
\label{eqn:fpv}
\end{align}
where we start the iteration with $u^{(0)}(\bx) \equiv v^{(0)}(\bx) \equiv 0$. Note 
to assure convergence of the fixed point iteration, it is crucial to 
include all linear terms in the definition of $H_k$ and $H_{2k}$,
in particular $4\pi\eta_0(\bx) + 4\pi\eta^{(1)}(\bx)$.

The presence of inverses in (\ref{eqn:fpu})--(\ref{eqn:fpv}) means that we have to 
solve the PDEs with the operators (\ref{eqn:hk})--(\ref{eqn:h2k}) in the 
whole space, both inside and outside the rectangular region $V$. 
This can be done by discretizing (\ref{eqn:hk})--(\ref{eqn:h2k}) 
inside $V$ and then placing a perfectly matched layer (PML) around it
to account for $\mathbb{R}^d \setminus V$.
To that effect we replace the operators $H_k$ and $H_{2k}$ in 
(\ref{eqn:fpu})--(\ref{eqn:fpv}) with their PML counterparts. 
Following \cite{chen2013optimal}, in the case $d=2$ we define the
PML analogues of $H_k$ and $H_{2k}$ as
\begin{align}
H_{k}^{\mbox{\tiny PML}} & = \frac{\partial}{\partial
  x}\left(\frac{e_y(\bx)}{e_x(\bx)} \frac{\partial}{\partial x}\right)
+ \frac{\partial}{\partial y}\left(\frac{e_x(\bx)}{e_y(\bx)}
\frac{\partial}{\partial y}\right) + k^2 e_x(\bx) e_y(\bx) ( 1 +
4\pi\eta_0(\bx) + 4\pi\eta^{(1)}(\bx) ),
\label{eqn:hkpml} \\
H_{2k}^{\mbox{\tiny PML}} & = \frac{\partial}{\partial
  x}\left(\frac{e_y(\bx)}{e_x(\bx)} \frac{\partial}{\partial x}\right)
+ \frac{\partial}{\partial y}\left(\frac{e_x(\bx)}{e_y(\bx)}
\frac{\partial}{\partial y}\right) + (2k)^2 e_x(\bx) e_y(\bx) ( 1 +
4\pi\eta_0(\bx) + 4\pi\eta^{(1)}(\bx) ),
\label{eqn:h2kpml}
\end{align}
where $e_x(\bx)$ and $e_y(\bx)$ are defined by
\begin{equation}
e_x(\bx) = \left\{\begin{tabular}{ll} $1 - i a_0 \left(
\dfrac{d_x(\bx)}{L_x} \right)^2$, & if $\bx \in V_x^{\mbox{\tiny
    PML}}$ \\ $1$, & otherwise
\end{tabular}\right., 
\; e_y(\bx) = \left\{\begin{tabular}{ll} $1 - i a_0 \left(
\dfrac{d_y(\bx)}{L_y} \right)^2$, & if $\bx \in V_y^{\mbox{\tiny
    PML}}$ \\ $1$, & otherwise
\end{tabular}\right.
\end{equation}
Here $V_x^{\mbox{\tiny PML}}$ and $V_y^{\mbox{\tiny PML}}$ each contain 
two PML layers that we surround $V$ with in $x$ and $y$ directions
respectively. The layers in $V_x^{\mbox{\tiny PML}}$ have widths $L_x$,
the layers in $V_y^{\mbox{\tiny PML}}$ have widths $L_y$. In the numerical
experiments we take $L_x = L_y = 1.5 \lambda$. The functions 
$d_x(\bx)$ and $d_y(\bx)$ compute the distances from a point $\bx$
in the corresponding PML layer to $\partial V$. The constant 
$a_0 = 1.79$ is chosen according to \cite{chen2013optimal}.

Finally, to apply the fixed point iteration 
(\ref{eqn:fpu})--(\ref{eqn:fpv}) numerically, we discretize 
(\ref{eqn:hkpml})--(\ref{eqn:h2kpml}) on a tensor product
finite difference grid with a five-point stencil. Both 
$H_{k}^{\mbox{\tiny PML}}$ and $H_{2k}^{\mbox{\tiny PML}}$
are discretized on the same grid. Thus, the grid should be
refined enough to properly resolve the higher wavenumber 
operator $H_{2k}^{\mbox{\tiny PML}}$. In the numerical 
experiments we take equal grid steps in the $x$ and $y$ 
directions $h_x=h_y={\lambda}/{20}$.

\end{document}